\documentclass{amsart}
\usepackage{amscd}

\def\nonsep{non-separating blow-up sequence}
\def\toric{separating blow-up sequence}
\def\why{} \def\regular{regular (i.e., generic)\def\regular{regular}}
\def\regular{generic\def\why{ (since generic fibres are also called
    ``regular'')}} 
\usepackage{xypic} \xyoption{all}
\usepackage{graphicx} \newcommand{\bbC}{{\mathbb C}}
\newcommand{\bbP}{{\mathbb P}} \newcommand{\bbZ}{{\mathbb Z}}
\newcommand{\C}{{\mathcal C}} \newcommand{\isom}{\cong}
\newtheorem{lemma}{Lemma}[section]
\newtheorem{prop}[lemma]{Proposition}
\newtheorem{thm}{Theorem}[section]
\newtheorem{other}{Theorem}[section] 
\newtheorem{cor}[other]{Corollary} 

\begin{document}
\title{Rational polynomials of simple type}
\author{Walter D. Neumann}
\address{Department of Mathematics\\Barnard College\\Columbia University\\
NY 10027\\USA}
\email{neumann@math.columbia.edu}
\author{Paul Norbury}
\address{Department of Mathematics and Statistics\\University of
Melbourne\\Australia 3010}
\email{norbs@ms.unimelb.edu.au}
\thanks{This research was supported by the Australian Research Council}

\begin{abstract}
We classify two-variable polynomials which are rational of simple
type. These are precisely the two-variable polynomials with trivial
homological monodromy.
%  We classify rational polynomials of simple type in two variables
%  whose generic fibres do not form an isotrivial family, and correct
%  the isotrivial case, which had been done earlier. This completes the
%  classification of all rational polynomials of simple type in two
%  variables.
\end{abstract}

\maketitle

\section{Introduction}
A polynomial map $f\colon\bbC^2\to \bbC$ is \emph{rational} if its
\regular{} fibre, and hence every fibre, is of genus
zero. It is of \emph{simple type} if, when extended to a morphism
$\tilde f\colon X\to\bbP^1$ of a compactification $X$ of $\bbC^2$, the
restriction of $\tilde f$ to each curve $C$ of the compactification
divisor $D=X-\bbC^2$ is either degree $0$ or $1$. The curves $C$ on
which $\tilde f$ is non-constant are called \emph{horizontal curves},
so one says briefly ``each horizontal curve is degree 1''.

The classification of rational polynomials of simple type gained some
new interest through the result of Cassou-Nogues, Artal-Bartolo, and
Dimca \cite{artal-cassou-dimca} that they are precisely the
polynomials whose homological monodromy is trivial (it suffices that
the homological monodromy at infinity be trivial by an observation of
Dimca).  

A classification appeared in \cite{MSuGen}, but it is incomplete.  It
implicitly assumes trivial \emph{geometric} monodromy (on page 346,
lines 10--11). Trivial geometric monodromy implies isotriviality
(\regular{} fibres pairwise isomorphic) and turns out to be equivalent
to it for rational polynomials of simple type.  The classification in
the non-isotrivial case was announced in the final section of
\cite{NNoMo}.  The main purpose of this paper is to prove it.  But we
recently discovered that there are also isotrivial rational
polynomials that are not in \cite{MSuGen}, so we have added a
classification for the isotrivial case using our methods.  This case
can also be derived from Kaliman's classification \cite{KalPol} of
\emph{all} isotrivial polynomials.  The fact that his list includes
rational polynomials of simple type that are not in \cite{MSuGen}
appears not to have been noticed before (it also includes rational
polynomials not of simple type).

In general, the classification of polynomial maps
$f:\bbC^2\rightarrow\bbC$ is an open problem with extremely rich
structure.  One notable result is the theorem of Abhyankar-Moh and
Suzuki \cite{AMoEmb,SuzPro} which classifies all polynomials with one
fibre isomorphic to $\bbC$.  The analogous result for the next
simplest case, where one fibre is isomorphic to $\bbC^*$, is open
except in special cases when the genus of the \regular{} fibre of the
polynomial is given.  Kaliman \cite{KalRat} classifies all {rational}
polynomials with one fibre isomorphic to $\bbC^*$.  

The basic tool we use in our study of rational polynomials is to
associate to any rational polynomial $f:\bbC^2\rightarrow\bbC$ a
compactification $X$ of $\bbC^2$ on which $f$ extends to a
well-defined map $\tilde{f}:X\rightarrow\bbP^1$ together with a map
$X\rightarrow\bbP^1\times\bbP^1$.  The map to $\bbP^1\times\bbP^1$ is
not in general canonical.  We will exploit the fact that for a
particular class of rational polynomials, there is an almost canonical
choice.

Although we give explicit polynomials, the classification is initially
presented in terms of the splice diagram for the link at infinity of a
\regular{} fibre of the polynomial (Theorem \ref{th:main}).  This is
called the \emph{regular splice diagram} for the polynomial\why.  See
\cite{NeuCom} for a description of the link at infinity and its splice
diagram.  The regular splice diagram determines the embedded topology
of a \regular{} fibre and the degree of each horizontal curve.  Hence
we can speak of a ``rational splice diagram of simple type''.

The first author has asked if the moduli space of polynomials with
given regular splice diagram is connected.  For a rational splice
diagram of simple type we find the answer is ``yes''.  We describe the
moduli space for our polynomials in Theorem \ref{th:defspace} and use
it to help give explicit normal forms for the polynomials.  We also
describe how the topology of the irregular fibres varies over the
moduli space.

The more general problem of classifying all rational polynomials,
which would cover much of the work mentioned above, is still an open
and interesting problem.  It is closely related to the problem of
classifying birational morphisms of the complex plane since a
polynomial is rational if and only if it is one coordinate of a
birational map of the complex plane.  Russell \cite{RusGoo} calls this
a ``field generator'' and defines a good field generator to be a
rational polynomial that is one coordinate of a birational morphism of
the complex plane.  A rational polynomial is good precisely when its
resolution has at least one degree one horizontal curve,
\cite{RusGoo}.  Daigle \cite{DaiBir} studies birational morphisms
$\bbC^2\rightarrow\bbC^2$ by associating to a compactification $X$ of
the domain plane a canonical map $X\rightarrow\bbP^2$.  A birational
morphism is then given by a set of curves and points in $\bbP^2$
indicating where the map is not one-to-one.  The approach we use in
this paper is similar.

The full list of rational polynomials $f\colon\bbC^2\to\bbC$ of simple
type is as follows. We list them up to polynomial automorphisms of
domain $\bbC^2$ and range $\bbC$ (so-called ``right-left
equivalence'').

\begin{thm}\label{th:summary}
  Up to right-left equivalence a rational polynomial $f(x,y)$ of
  simple type has one of the following forms $f_i(x,y)$, $i=1$, $2$,
  or $3$.
\begin{align*}
f_1(x,y)=&x^{q_1}s^q+x^{p_1}s^p
\prod_{i=1}^{r-1}(\beta_i-x^{q_1}s^q)^{a_i}&(r\ge2)\phantom{.}\\
f_2(x,y)=&x^{p_1}s^p
\prod_{i=1}^{r-1}(\beta_i-x^{q_1}s^q)^{a_i}&(r\ge1)\phantom{.}\\
f_3(x,y)=&y\prod_{i=1}^{r-1}(\beta_i-x)^{a_i}+h(x)&(r\ge1).
\end{align*}
Here:

$0\le q_1<q,\quad 0\le p_1<p,\quad \left|
  \begin{matrix}
    p&p_1\\q&q_1
  \end{matrix}\right|=\pm1 $;

$s=yx^k+P(x),\text{ with $k\ge1$ and $P(x)$ a polynomial of
  degree $<k$}$;

$a_1,\dots,a_{r-1}\text{ are positive integers}$;

$\beta_1,\dots,\beta_{r-1}\text{ are distinct elements of }\bbC^*$;

$h(x)\text{ is a polynomial of degree }< \sum_1^{r-1}a_i$.

Moreover, if $g_1(x,y)=g_2(x,y)=x^{q_1}s^q$ and $g_3(x,y)=x$ then
$(f_i,g_i)\colon\bbC^2\to\bbC^2$ is a birational morphism for
$i=1,2,3$. In fact, $g_i$ maps a \regular{} fibre $f_i^{-1}(t)$
biholomorphically to $\bbC-\{0,t,\beta_1,\dots,\beta_{r-1}\}$,
$\bbC-\{0,\beta_1,\dots,\beta_{r-1}\}$, or
$\bbC-\{\beta_1,\dots,\beta_{r-1}\}$, according as $i=1,2,3$. Thus
$f_1$ is not isotrivial and $f_2$ and $f_3$ are.
\end{thm}

In \cite{MSuGen} the isotrivial case is subdivided into seven
subcases, but these do not include any $f_2(x,y)$ with $p,q,p_1,q_1$
all $>1$.

\section{Resolution}
Given a polynomial $f:\bbC^2\rightarrow\bbC$, extend it to a map
$\bar{f}:\bbP^2\rightarrow\bbP^1$ and resolve the points of
indeterminacy to get a regular map $\tilde{f}:X\rightarrow\bbP^1$ that
coincides with $f$ on $\bbC^2\subset X$.  We call $D=X-\bbC^2$ the
divisor at infinity.  The divisor $D$ consists of a connected union of
rational curves.  An irreducible component $E$ of $D$ is {\em
  horizontal} if the restriction of $\tilde{f}$ to $E$ is not a
constant mapping.  The \emph{degree of a horizontal curve} $E$ is the
degree of the restriction $\tilde{f}|E$.  Although the
compactification defined above is not unique, the horizontal curves
are essentially independent of choice.

Note that a \regular{} fibre $F_c:=f^{-1}(c)$ is a punctured Riemann
surface with punctures precisely where $\overline F_c$ meets a
horizontal curve. Thus $f$ has simple type if and only if $\overline
F_c$ meets each horizontal curve exactly once, so the number of
punctures equals the number of horizontal curves. For non-simple type
the number of punctures will exceed the number of horizontal curves.

We say that a rational polynomial is {\em ample} if it has at least
three degree one horizontal curves.  Those polynomials with no degree
one horizontal curves, or bad field generators \cite{RusGoo}, are
examples of polynomials that are not ample.  The classification of
Kaliman \cite{KalRat} mentioned in the introduction gives examples of
polynomials with exactly one degree one horizontal curve so they are
also not ample.  Nevertheless, ample rational polynomials will be the
focus of our study in this paper.  
We will classify all ample rational polynomials that are also of
simple type. 

\section{Curves in $\bbP^1\times\bbP^1$.}  \label{sec:cur}
If $\tilde{f}:X\rightarrow\bbP^1$ is a regular map with rational
fibres then $X$ can be blown down to a Hirzebruch surface, $S$, so
that $\tilde{f}$ is given by the composition of the sequence of
blow-downs $X\rightarrow S$ with the natural map $S\rightarrow\bbP^1$;
see \cite{BavdV} %\cite{GHaPri}, p.513 
for details.  Moreover, by first replacing
$X$ by a blown-up version of $X$ if necessary, we may assume that
$S=\bbP^1\times\bbP^1$
and the natural map to $\bbP^1$ is projection onto the first factor.

A rational polynomial $f\colon\bbC^2\to\bbC$, once compactified to
$\tilde f\colon X=\bbC^2\cup D\to \bbP^1$, may thus be given by
$\bbP^1\times\bbP^1$ together with instructions how to blow up
$\bbP^1\times\bbP^1$ to get $X$ and how to determine $D$ in $X$. For
this we give the following data:
\begin{itemize}
\item a collection $\mathcal C$ of irreducible rational curves in
  $\bbP^1\times\bbP^1$ including $L_\infty:=\infty\times\bbP^1$;
\item a set of instructions on how to blow up $\bbP^1\times\bbP^1$ to
  obtain $X$;
\item a sub-collection $\mathcal E$ of the curves of the exceptional
  divisor of $X\to\bbP^1\times\bbP^1$;
\end{itemize}
satisfying the condition:
\begin{itemize}
\item If $D$ is the union of the curves of $\mathcal E$ and the proper
  transforms of the curves of $\mathcal C$ then $X-D\isom\bbC^2$;
\end{itemize}

If $C\subset\bbP^1\times\bbP^1$ is an irreducible algebraic curve we
associate to it the pair of integers $(m,n)$ given by degrees of the
two projections of $C$ to the factors of $\bbP^1\times\bbP^1$.
Equivalently, $(m,n)$ is the homology class of $C$ in terms of
$H_2(\bbP^1\times\bbP^1)=\bbZ\oplus\bbZ$.  We call $C$ an $(m,n)$
curve.  The intersection number of an $(m,n)$ curve $C$ and an
$(m',n')$ curve $C'$ is $C\cdot C'=mn'+nm'$.

The above collection $\mathcal C$ of curves in $\bbP^1\times\bbP^1$ will
consist of some \emph{vertical curves} (that is, $(0,1)$ curves; one
of these is $L_\infty$) and some other curves. These non-vertical
curves give the horizontal curves for $f$, so they all have $m=1$ if
$f$ is of simple type. Note that a $(1,n)$ curve is necessarily smooth
and rational (since it is the graph of a morphism $\bbP^1\to\bbP^1$).

The image in $\bbP^1\times\bbP^1$ of the fibre over infinity is the
$(0,1)$ curve $L_{\infty}$ and the image of a degree $m$ horizontal
curve is an $(m,n)$ curve.  This view allows one to see as follows a
geometric proof of the result of Russell \cite{RusGoo} that a rational
polynomial $f$ is good precisely when its resolution has at least one
degree one horizontal curve.  A degree one horizontal curve
for $f$ has image in $\bbP^1\times\bbP^1$ given by a $(1,n)$ curve.
Call this image $C$ and let $P$ be its intersection with $L_\infty$.
The $(1,n)$ curves that do not intersect $C-P$ form a $\bbC$--family
that sweeps out $\bbP^1\times\bbP^1- (L_\infty\cup C)$ so they
lead to a map $X\to\bbP^1$ which takes values in $\bbC$ at points that
do not lie over $L_\infty\cup C$. Restricting to $\bbC^2=X-D$ we
obtain a meromorphic function $g_1$ that has poles only at points that
belong to exceptional curves that were blown up on $C$ (and do not
belong to $\mathcal E$). However the polynomial $f$ is constant on
each such curve, so if $c_1, \dots, c_k$ are the values that $f$ takes
on these curves, then $g:=g_1(f-c_1)^{a_1}\dots(f-c_k)^{a_k}$ will
have no poles, and hence be polynomial, for $a_1, \dots, a_k$
sufficiently large. Then $(f,g)$ is the desired birational morphism
$\bbC^2\to \bbC^2$. For the converse, given a birational morphism
$(f,g)\colon\bbC^2\to\bbC^2$, we compactify it to a morphism $(\tilde
f,\tilde g)\colon X\to\bbP^1\times\bbP^1$.  Then the proper transform
of $\bbP^1\times\infty$ is the desired degree one horizontal curve for
$f$.

We shall use the usual encoding of the topology of $D$ by the dual
graph, which has a vertex for each component of $D$, an edge when two
components intersect, and vertex weights given by self-intersection
numbers of the components of $D$. We will sometimes speak of the
\emph{valency} of a component $C$ of $D$ to mean the valency of the
corresponding vertex of the dual graph, that is, the number of other
components that $C$ meets.

The approach we will take to get rational polynomials will be to start
with any collection $\C$ of $k$ curves in $\bbP^1\times\bbP^1$ and see
if we can produce a divisor at infinity $D$ for a map from $\bbC^2$ to
$\bbC$.  In order to get a divisor at infinity we must blow up
$\bbP^1\times\bbP^1$, say $m$ times, and include some of the resulting
exceptional curves in the collection so that this new collection gives
a divisor $D$ whose complement is $\bbC^2$. The exceptional curves
that we ``leave behind'' (i.e., do not include in $D$) will be called
\emph{cutting divisors}.

\begin{lemma}  \label{th:propD}
(i) $D$ must have $m+2$ irreducible components, so we must
include $m-k+2$ of the exceptional divisors in the collection leaving
$k-2$ behind as cutting divisors;

(ii) $D$ must be connected and have no cycles;

(iii) $D$ must reduce to one of the ``Morrow configurations'' by a
sequence of blow-downs. The Morrow configurations are the
configurations of rational curves with dual graphs of one of the
following three types,
in which, in the last case, after replacing the central $(n,0,-n-1)$ by
a single $(-1)$ vertex the result should blow down to a single
$(+1)$ vertex by a sequence of blow-downs:
$$
 \xymatrix@R=6pt@C=30pt@M=0pt@W=0pt@H=0pt{
_{1}\\
\circ}$$ $$
\xymatrix@R=6pt@C=30pt@M=0pt@W=0pt@H=0pt{
_{0}&_{l}\\
\circ\ar@{-}[r]&\circ}$$ $$
\xymatrix@R=6pt@C=30pt@M=0pt@W=0pt@H=0pt{
_{l_m}&_{\cdots}&_{l_1}&_{n}&_{0}&_{-n-1}&_{t_1}&_{\cdots}&_{t_k}\\
\circ\ar@{.}[rr]&&\circ\ar@{-}[r]&\circ\ar@{-}[r]&\circ\ar@{-}[r]&
\circ\ar@{-}[r]&\circ\ar@{.}[rr]&&\circ}
$$

\vspace{6pt}These conditions are also sufficient that $X-D\isom \bbC^2$.
\end{lemma}
\begin{proof}
The first property follows from the fact that each blow-up
increases the rank of second homology by $1$. Thus $H_2(X)$ has rank
$m+2$, so $D$ must have $m+2$ irreducible components.
Notice that this implies easily the well-known result
\cite{KalTwo,MSuGen,SuzPro} that
\[\delta-1=\sum_{a\in\bbC}(r_a-1),\]
where $\delta$ is the number of horizontal curves of $f$ and $r_a$ is
the number of irreducible components of $f^{-1}(a)$.  (Both sides are
equal to $k-1-\{$number of finite curves at infinity$\}$.)

The second property follows from the third property.  For the third
property and sufficiency see \cite{Morrow, Rama}.
\end{proof}

Now assume that $\tilde{f}$ has at least three degree one horizontal
curves.  Take these three horizontal curves and use them to map $X$ to
$\bbP^1\times\bbP^1$ as follows.  The three horizontal curves define
three points in a \regular{} fibre of $\tilde{f}$. We can map this
\regular{} fibre to $\bbP^1$ by mapping these three points to
$0,1,\infty\in\bbP^1$. This defines a map from a Zariski open set of
$X$ to $\bbP^1$ which then extends to a map $\pi$ from $X$ to
$\bbP^1$.  If $\pi$ is not a morphism then we blow up $X$ to get a
morphism.  Rather than introducing further notation for this blow-up
we will assume we began with this blow-up and call it $X$.  Together
with the map $\tilde{f}$ this gives us the desired morphism
\[X\stackrel{(\tilde{f},\pi)}{\longrightarrow}\bbP^1\times\bbP^1\]
with the property that the three horizontal curves map to $(1,0)$
curves.  

If all horizontal curves for $f$ are of type $(1,0)$ then the \regular{}
fibres form an isotrivial family (briefly ``$f$ is isotrivial'').
Thus if $f$ is of simple type but not isotrivial, there must be a
horizontal curve of type $(1,n)$ in $\C$ with $n>0$.  From now on,
therefore, we assume that there are at least three $(1,0)$ curves and
at least one $(1,n)$ curve in $\C$ with $n>0$.

\begin{lemma}\label{le:beyond}
  Any curve of $D$ that is beyond a horizontal curve from the point of
  view of $\tilde L_\infty$ has self-intersection number $\le-2$.
% or  self-intersection number $-1$ and it meets at least three other
%  curves of $D$.
\end{lemma}
\begin{proof}
  If the curve is an exceptional curve then it has self-intersection
  $\le-1$. If $-1$, then the curve must have valency at least three
  (since any $-1$ exceptional curve that could be blown down is a
  cutting divisor). Any three adjacent curves must include two
  horizontal curves, which contradicts the fact that the dual graph of
  $D$ has no cycles. If the curve is not exceptional then it is the
  proper transform of a vertical curve. But we must have blown up at
  least three times on the vertical curve to get rid of cycles in the
  dual graph of $D$ so in this case the self-intersection is $\le-3$.
\end{proof}

\subsection{Horizontal curves}
The next few lemmas will be devoted to finding restrictions on the
horizontal curves in the configuration $\C\subset\bbP^1\times\bbP^1$,
culminating in Proposition~\ref{th:canon}.

\begin{lemma}   \label{th:n=1}
  A horizontal curve of type $(1,n)$ in $\C$ must be of type $(1,1)$.
\end{lemma}
\begin{proof}  
  Assume we have a horizontal curve $C\in \C$ of type $(1,n)$ with
  $n>1$.  It intersects each of the three $(1,0)$ curves $n$ times
  (counting with multiplicity) so in order to break
  cycles---Lemma~\ref{th:propD} ~(ii)---we have to blow up at least
  $n$ times on each $(1,0)$ horizontal curve, so the proper transforms
  of the three $(1,0)$ curves have self-intersection at most $-n$ and
  the proper transform of the $(1,n)$ curve has self-intersection at
  most $2n-3n=-n$.
  
  By Lemma~\ref{th:propD} ~(iii), $D$ must reduce to a Morrow
  configuration by a sequence of blow-downs.  Thus $D$ must contain a
  $-1$ curve $E$ that blows down.  
  By Lemma
  \ref{le:beyond}, the curve $E$ must be a proper transform of a
  horizontal curve.  The proper transform of each $(1,0)$ curve has
  self-intersection at most $-n<-1$.  Thus $E$ must come from one of
  the $(1,*)$ horizontal curves.  As mentioned above, the proper
  transform of a $(1,k)$ curve has self-intersection $\leq -k$ so $E$
  must be the proper transform of a $(1,1)$ curve, $E_0$.  But $E_0$
  would intersect $C$, the $(1,n)$ curve, $2n$ times and hence
  $E.E\leq 2-2n<-1$ since $n>1$.  This is a contradiction so any
  horizontal curve of type $(1,n)$ must be a $(1,1)$ curve.
\end{proof}

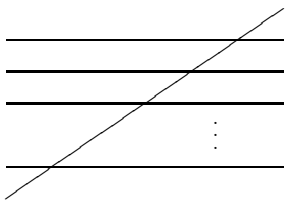
\begin{figure}[ht]
\def\Vdots{\raise5pt\hbox{$\vdots$}}
$$\xymatrix@R=12pt@C=15pt@M=0pt@W=0pt@H=0pt{
\ar@{-}&&&&&&&\\
\ar@{-}[rrrrrrr]&&&&&&&\\
\ar@{-}[rrrrrrr]&&&&&&&\\
\ar@{-}[rrrrrrr]&&&&&\ar@{}[dd]^(.3){.}^{.}^(.7){.}&&\\
\ar@{-}&&&&&&&\\
\ar@{-}[rrrrrrr]&&&&&&&\\
\ar@{-}[rrrrrrruuuuuu]}$$
\caption{Configuration of horizontal curves.}
\label{fig:conf1}
\end{figure}

Hence, the horizontal curves consist of a collection of $(1,0)$ curves
and $(1,1)$ curves.  Figure~\ref{fig:conf1} shows an example of a
possible configuration of horizontal curves in $\bbP^1\times\bbP^1$.

\begin{lemma}   \label{th:triple}
  $\tilde{L}_{\infty}\cdot \tilde{L}_{\infty}=-1$.
\end{lemma}
\begin{proof}  We blow up at a point on $L_{\infty}$
  precisely when at least two horizontal curves meet in a common point
  there.  In general, if a horizontal curve meets $L_{\infty}$ with a
  high degree of tangency then we blow up repeatedly there.  But,
  since all horizontal curves are $(1,0)$ and $(1,1)$ curves, they
  meet $L_{\infty}$ transversally, so a point on $L_{\infty}$ will be
  blown up at most once.
  
  If there are two such points to be blown up, then after blowing up
  there will be (in the dual graph) two non-neighbouring $-1$ curves
  with valency $>2$.  The complement of such a configuration cannot be
  $\bbC^2$.  This is proven by Kaliman \cite{KalTwo} as Corollary 3.
  Actually the result is stated for two $-1$ curves of valency 3 but
  it applies to valency $\ge3$.
  
  Thus, at most one point on $L_{\infty}$ is blown up and
  $\tilde{L}_{\infty}\cdot \tilde{L}_{\infty}=0$ or $-1$. We must show
  $0$ cannot occur.
  
  Since there are at least four horizontal curves, if
  $\tilde{L}_{\infty}\cdot \tilde{L}_{\infty}=0$, then
  $\tilde{L}_{\infty}$ has valency at least $4$ and every other curve
  has negative self-intersection.  Furthermore, the only possible $-1$
  curves must be horizontal curves, and these intersect
  $\tilde{L}_{\infty}$ in $D$.  As we attempt to blow down $D$ to get
  to a Morrow configuration, the only curves that can be blown down
  will always be adjacent to $\tilde L_{\infty}$.  Thus the
  intersection number of $\tilde L_\infty$ will become positive and
  all other intersection numbers remain negative, so a Morrow
  configuration cannot be reached.  Hence, $\tilde{L}_{\infty}\cdot
  \tilde{L}_{\infty}=-1$.
\end{proof}

\begin{figure}
\begin{picture}(200,100)
\put(50,40){\framebox(20,20){$A$}}
\put(100,70){\framebox(20,20){$B_1$}}
\put(100,10){\framebox(20,20){$B_2$}}
\put(110,50){\circle{6}e}
\put(70,50){\line(1,0){37}}
\put(110,70){\line(0,-1){17}}
\put(110,30){\line(0,1){17}}
\end{picture}
\caption{The branches $B_1$ and $B_2$ consist of curves of 
self-intersection $<-1$ and $e\geq -1$.}
\label{fig:neg}
\end{figure}
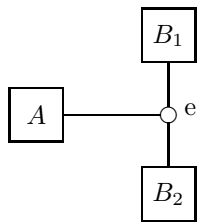

\begin{lemma}  \label{th:neg}
  A configuration of curves that contains two branches consisting of
  curves of self-intersection $<-1$ that meet at a valency $>2$ curve
  of self-intersection greater than or equal to $-1$ as in
  Figure~\ref{fig:neg} (where the meeting curve is drawn with valency
  3 for convenience) cannot be blown down to a Morrow configuration.
\end{lemma}
\begin{proof}
  Since the two branches consist of curves of self-intersection $<-1$,
  they cannot be reduced before the other branches are reduced.  If
  the rest of the configuration of curves is blown down first then the
  valency $>2$ curve becomes a valency $2$ curve with non-negative
  self-intersection and no more blow-downs can be done. Since there is
  no $0$ curve, we have not reached a Morrow configuration.
\end{proof}

\begin{lemma}   \label{th:inunion}
  The intersection of any two $(1,1)$ curves in $\C$ consists of two
  distinct points contained in the union of the $(1,0)$ curves in
  $\C$.
\end{lemma}
\begin{proof} 
  We will assume otherwise and reduce to the situation of
  Lemma~\ref{th:neg} to give a contradiction.  Thus, assume that two
  $(1,1)$ curves do not intersect in two points contained in the union
  of the $(1,0)$ curves.  Then in order to break cycles these curves
  must be blown up at least four times---once each for at least three
  of the $(1,0)$ curves and at least another time for the intersection
  of the two $(1,1)$ curves.  Thus they have self-intersection $<-1$.
  
  {\em Case 1}: Suppose two $(1,1)$ curves meet on $L_{\infty}$.  Then
  after blowing up (twice if the $(1,1)$ curves meet at a tangent),
  the exceptional curves are retained and the final exceptional curve
  has self-intersection $-1$, valency 3 and two branches, which we
  will call $B_1$ and $B_2$, consisting of the proper transforms of
  the two $(1,1)$ curves and any other curves beyond these proper
  transforms all of which have self-intersection $<-1$.  Thus we are
  in the situation of Lemma~\ref{th:neg} and we get a contradiction.
  
  {\em Case 2}: Suppose two $(1,1)$ curves meet $L_{\infty}$ at
  distinct points.  Then at least one of the $(1,1)$ curves, $D$, must
  meet $L_{\infty}$ at a point away from the $(1,0)$ curves by
  Lemma~\ref{th:triple}.  Also one of the $(1,0)$ curves, $H$, must
  meet $L_{\infty}$ away from the $(1,1)$ curves and contain at least
  two points where it intersects the $(1,1)$ curves and thus have
  self-intersection $<-1$ after blowing up to break cycles.  We are
  once more at the situation of Lemma~\ref{th:neg} where the valency
  $>2$ curve is $\tilde{L}_{\infty}$ which has self-intersection $-1$
  by Lemma~\ref{th:triple}, and the branches $B_1$ and $B_2$ are the
  proper transform of $D$ and any curves beyond it, respectively the
  proper transform of $H$ and any curves beyond it.  Thus we have a
  contradiction.
  
  Notice that both cases apply to two $(1,1)$ curves that may
  intersect at a tangent point, and shows that this situation is
  impossible.
\end{proof}

\begin{lemma}   \label{th:three}
  If there is more than one $(1,1)$ curve in $\C$ then there are
  exactly three $(1,0)$ horizontal curves in $\C$.
\end{lemma}
\begin{proof} 
  Assume that there are more than three $(1,0)$ horizontal curves in
  $\C$ and at least two $(1,1)$ curves, say $C_1$ and $C_2$.
  
  {\em Case 1}: $C_1$ and $C_2$ meet on $\tilde L_\infty$. Then they
  meet each of at least two $(1,0)$ curves in distinct points, so
  after blowing up to destroy cycles, these $(1,0)$ curves have
  self-intersection number $\le-2$ and Lemma \ref{th:neg} applies.
 
  {\em Case 2}: $C_1$ and $C_2$ meet $\tilde L_\infty$ at distinct
  points. Then one of them, say $C_1$, meets $\tilde L_\infty$ at a
  point not on a $(1,0)$ curve by Lemma \ref{th:triple}. At least one
  $(1,0)$ curve $C_3$ meets $C_1$ and $C_2$ in distinct points. After
  breaking cycles, $C_1$ and $C_3$ have self-intersections $\le -2$ so
  Lemma \ref{th:neg} applies again.
\end{proof}

\begin{lemma}   \label{th:common}
  A family of $(1,1)$ horizontal curves in $\C$ must pass through a
  common pair of points.
\end{lemma}
\begin{proof}  
  The statement is trivial for one $(1,1)$ horizontal curve so assume
  there are at least two $(1,1)$ horizontal curves in $\C$.  By the
  previous lemma, there are exactly three $(1,0)$ horizontal curves.
  
  If there are exactly two $(1,1)$ horizontal curves in $\C$ then the
  lemma is clear since the curves cannot be tangent by
  Lemma~\ref{th:inunion}.
  
  When there are more than two $(1,1)$ curves in $\C$, apply
  Lemma~\ref{th:inunion} to two of them.  If another $(1,1)$
  horizontal curve in $\C$ does not intersect these two $(1,1)$ curves
  at their common two points of intersection then, by
  Lemma~\ref{th:inunion}, it must meet both these $(1,1)$ curves at
  the third $(1,0)$ horizontal curve of $\C$. So the first two $(1,1)$
  curves would meet there, which is a contradiction.
\end{proof}

\begin{prop}   \label{th:canon}
  Any configuration of horizontal curves in $\C$ is equivalent to one
  of the form in Figure~\ref{fig:conf1}.
\end{prop}
\begin{proof}  
  By assumption and Lemma~\ref{th:n=1} there are at least three
  $(1,0)$ horizontal curves and some $(1,1)$ horizontal curves in
  $\C$.  If there is exactly one $(1,1)$ horizontal curve then the
  proposition is clear.  If there is more than one $(1,1)$ horizontal
  curve, then by Lemmas~\ref{th:three} and \ref{th:common} there are
  precisely three $(1,0)$ horizontal curves and two of the $(1,0)$
  horizontal curves contain the common intersection of the $(1,1)$
  curves.  Each $(1,1)$ curve also contains a distinguished point
  where the curve meets the third $(1,0)$ horizontal curve.  A Cremona
  transformation can bring such a configuration to that in
  Figure~\ref{fig:conf1} by blowing up at the two points of
  intersection of the $(1,1)$ curves and blowing down the two vertical
  lines containing the two points.  This sends two of the $(1,0)$
  horizontal curves and each $(1,1)$ curve to $(1,0)$ horizontal
  curves and one of the $(1,0)$ curves to a $(1,1)$ curve that
  intersects each of the other horizontal curves exactly once.  Note
  that since we blow up $\bbP^1\times\bbP^1$ to get the polynomial
  map, two configurations of curves $\C,\C'$ in $\bbP^1\times\bbP^1$
  related by a Cremona transformation give rise to the same
  polynomial, so we are done.
\end{proof}

\subsection{The configuration $\C$}
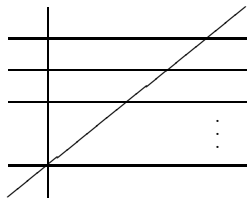
\begin{figure}[ht]
$$\xymatrix@R=12pt@C=15pt@M=0pt@W=0pt@H=0pt
{\ar@{-}&\ar@{-}[dddddd]&&&&&&\\
\ar@{-}[rrrrrr]&&&&&&&\\
\ar@{-}[rrrrrr]&&&&&&&\\
\ar@{-}[rrrrrr]&&&&&\ar@{}[dd]^(.3){.}^{.}^(.7){.}&&\\
\ar@{-}&&&&&&&\\
\ar@{-}[rrrrrr]&&&&&&&\\
\ar@{-}[rrrrrruuuuuu]&&&&&&&}$$
\caption{Configuration $\C$ with $r+2$ horizontal curves.}
\label{fig:confC}
\end{figure}

The image $\C$ of $D\subset X\rightarrow\bbP^1\times\bbP^1$ will
consist of the configuration of horizontal curves in
Figure~\ref{fig:conf1} plus some $(0,1)$ vertical curves.  The next
two lemmas show that in fact the only $(0,1)$ vertical curve we need
to include in $\C$ is $L_{\infty}$ and furthermore that $\C$ can be
given by Figure~\ref{fig:confC}.

\begin{lemma}\label{le:list}
  The configuration $\C$ appears in Figure~\ref{fig:list} or
  Figure~\ref{fig:confC}.
\end{lemma}
\begin{proof}
  Let $r+2$ denote the number of horizontal curves and $k+1$ denote
  the number of $(0,1)$ vertical curves in $\C$.  Thus $\C$ consists
  of $k+r+3$ irreducible components and by Lemma~\ref{th:propD} ~(i),
  when blowing up to get $D$ from $\C$ we must leave $k+r+1$
  exceptional curves behind as cutting divisors.
  
  By Lemma~\ref{th:propD} ~(ii) we must break all cycles.  The minimum
  number of cutting divisors needed to do this is $kr+k+r-2{\rm
    min}\{k,r\}$.  This is because each of the $k$ $(0,1)$ vertical
  curves different from $L_{\infty}$ must be separated from all but
  one of the $r+1$ $(1,0)$ horizontal curves, so we need $kr$ cutting
  divisors.  Also, the $(1,1)$ horizontal curve meets each of the
  $r+1$ $(1,0)$ horizontal curves and each of the $k$ $(0,1)$ vertical
  curves once, so that requires $k+r$ cutting divisors (by
  Lemma~\ref{th:triple} the $(1,1)$ curve must meet $L_{\infty}$ at a
  triple point with a $(1,0)$ horizontal curve, so this intersection
  does not produce a cycle to be broken).  We would thus require
  $kr+k+r$ cutting divisors except that the $(1,1)$ curve may pass
  through intersections of the $(1,0)$ horizontal curves and the
  $(0,1)$ vertical curves, so some of the cutting divisors may
  coincide.  The most such intersections possible is ${\rm
    min}\{k,r\}$ and we have then over-counted required cutting
  divisors by $2{\rm min}\{k,r\}$. Hence we get at least $kr+k+r-2{\rm
    min}\{k,r\}$ cutting divisors.
  
  Since the number $k+r+1$ of cutting
  divisors is at least $kr+k+r-2{\rm min}\{k,r\}$, we have $k+r+1\geq
  kr+k+r-2{\rm min}\{k,r\}$, so
  \begin{equation}  \label{eq:ineq}
  1\geq k(r-2)\ \ {\rm and}\ \ 1\geq (k-2)r,\ \ \ k\geq 0,\ r\geq 2.
  \end{equation}
  The solutions of (\ref{eq:ineq}) are
  $(k,r)=\{(0,r),(1,2),(1,3),(2,2)\}$.
  
  Recall by Lemma~\ref{th:triple} that the $(1,1)$ curve must meet
  $L_{\infty}$ at a triple point with a $(1,0)$ horizontal curve.
  Furthermore, by keeping track of when either inequality in
  (\ref{eq:ineq}) is an equality, or one away from an equality, we can
  see that the $(1,1)$ curve must meet any other $(0,1)$ vertical
  curves at a triple point with a $(1,0)$ horizontal curve.  Thus, the
  only possible configurations for $\C$ are given in
  Figures~\ref{fig:list} and \ref{fig:confC}.
\end{proof}

\begin{figure}[ht]
$$\xymatrix@R=12pt@C=15pt@M=0pt@W=0pt@H=0pt{
\ar@{-}&\ar@{-}[dddddd]&&&&\ar@{-}[dddddd]&&&\ar@{-}&\ar@{-}[dddddd]&&&&\ar@{-}[dddddd]&&&\ar@{-}&\ar@{-}[dddddd]&&\ar@{-}[dddddd]&&\ar@{-}[dddddd]&&\\
\ar@{-}[rrrrrr]&&&&&&&&\ar@{-}[rrrrrr]&&&&&&&&\ar@{-}[rrrrrr]&&&&&&&\\
\ar@{-}&&&&&&&&\ar@{-}[rrrrrr]&&&&&&&&\ar@{-}&&&&&&&\\
\ar@{-}[rrrrrr]&&&&&&&&\ar@{-}[rrrrrr]&&&&&&&&\ar@{-}[rrrrrr]&&&&&&&\\
\ar@{-}&&&&&&&&\ar@{-}&&&&&&&&\ar@{-}&&&&&&&\\
\ar@{-}[rrrrrr]&&&&&&&&\ar@{-}[rrrrrr]&&&&&&&&\ar@{-}[rrrrrr]&&&&&&&\\
\ar@{-}[rrrrrruuuuuu]&&&&&&&&\ar@{-}[rrrrrruuuuuu]&&&&&&&&\ar@{-}[rrrrrruuuuuu]&&&&&&&}$$
\caption{Configuration $\C$.}
\label{fig:list}
\end{figure}

In the following lemmas we will exclude the configurations in
Figure~\ref{fig:list}.  Label the triple points in the first two
configurations of Figure~\ref{fig:list} by $P_{\infty}\in L_{\infty}$
and $P_1$, and in the third configuration by $P_\infty,P_1,P_2$.
Also, label the exceptional divisor obtained by blowing up the triple
point $P_i$ by $E_i$ and its proper transform by $\tilde{E}_i$.

\begin{lemma}   \label{th:confC}
If $E_i$ is a cutting divisor then the $(0,1)$ vertical curve
containing $P_i$ can be removed from $\C$ by a birational
transformation.
\end{lemma}
\begin{proof}
In each of the configurations of Figure~\ref{fig:list} we can perform
a Cremona transformation by blowing up $P_{\infty}$ and $P_i$ for
$i=1$ or $2$ and then blowing down $\tilde{L}_{\infty}$ and the proper
transform of the $(0,1)$ vertical curve that contains $P_i$.  The
exceptional divisors $E$ and $E_i$ become $(0,1)$ curves and the
$(0,1)$ vertical curve that contains $P_i$ becomes an exceptional
divisor in a new configuration $\C$.  When $E_i$ is a cutting divisor
this operation essentially removes a $(0,1)$ vertical curve from $\C$.
\end{proof}

\begin{lemma}
In a configuration from Figure~\ref{fig:list} with
$(k,r)\in\{(1,3),(2,2)\}$ at least one of the exceptional divisors $E_1$
or $E_2$ is a cutting divisor.
\end{lemma}
\begin{proof}
Suppose otherwise, that $E_1$ is not a cutting divisor and for
$(k,r)=(2,2)$ nor is $E_2$ a cutting divisor.  The exceptional curves
$E_i$ introduce an extra intersection and hence an extra cutting
divisor is required.  There is one such extra intersection in the
configuration with $(k,r)=(1,3)$ and two such extra intersections
in the configuration with $(k,r)=(2,2)$.  As mentioned in the
proof of Lemma~\ref{le:list} the solution $(k,r)=(1,3)$ gives
equality in (\ref{eq:ineq}) and so it cannot sustain an extra cutting
divisor.  Similarly the solution $(k,r)=(2,2)$ is $1$ away from
equality in (\ref{eq:ineq}) and so it cannot sustain two extra cutting
divisors.  Hence we get a contradiction and the lemma is proven.
\end{proof}

By the previous two lemmas we can simplify any configuration from
Figure~\ref{fig:list} to lie in Figure~\ref{fig:confC} or to be the
first configuration from Figure~\ref{fig:list} (the one with
$(k,r)=(1,2)$) with the requirement that $E_1$ is not a cutting
divisor.  It is this last case that we will now exclude.

The next three lemmas suppose that we have the first configuration
from Figure~\ref{fig:list} and that $E_1$ is not a cutting divisor.
We will denote the four horizontal curves by $H_i$, $i=1,\dots,4$, and
their proper transforms by $\tilde{H}_i$ where $H_4$ is the $(1,1)$
curve, $H_1$ contains $P_1$ and $H_3$ contains $P_{\infty}$.  Also
denote the $(1,0)$ vertical curve that contains $P_1$ by $L_1$ and its
proper transform by $\tilde{L}_1$.

\begin{lemma}    \label{th:atl}
  At least one of $\tilde{H}_1$ and $\tilde{H}_2$ and at least one of
  $\tilde{H}_3$ and $\tilde{H}_4$ has self-intersection $-1$.
\end{lemma}
\begin{proof}
  The proper transform of each horizontal curve has self-intersection
  less than or equal to $-1$ and all curves in $D$ beyond horizontal
  curves have self-intersection strictly less than $-1$.  If the two
  horizontal curves that meet $\tilde{L}_{\infty} $, $\tilde{H}_1$ and
  $\tilde{H}_2$, have self-intersection strictly less than $-1$, then
  since all curves beyond the two horizontal curves also have
  self-intersection strictly less than $-1$, and since
  $\tilde{L}_{\infty}$ has self-intersection $-1$ and valence $3$ this
  gives a contradiction by Lemma~\ref{th:neg}.  The same argument
  applies to $\tilde{H}_3$ and $\tilde{H}_4$ together with $E$.
\end{proof}

\begin{lemma}  \label{th:pair}
$\tilde{H}_4\cdot\tilde{H}_4= -1$ if and only if
$\tilde{H}_2\cdot\tilde{H}_2= -1$.
\end{lemma}
\begin{proof}
Since $L_1$ must be separated from at least one of $H_2$ and $H_3$
then at most one of $\tilde{H}_2\cdot\tilde{H}_2= -1$ and
$\tilde{H}_3\cdot\tilde{H}_3= -1$ can be true.  Similarly $E_1$ must
be separated from at least one of $H_1$ and $H_4$ so at most one of
$\tilde{H}_1\cdot\tilde{H}_1= -1$ and $\tilde{H}_4\cdot\tilde{H}_4=
-1$ can be true.  By Lemma~\ref{th:atl}, if
$\tilde{H}_2\cdot\tilde{H}_2\neq -1$ then
$\tilde{H}_1\cdot\tilde{H}_1= -1$ so $\tilde{H}_4\cdot\tilde{H}_4\neq
-1$.  Similarly, $\tilde{H}_1\cdot\tilde{H}_1\neq -1$ implies that
$\tilde{H}_2\cdot\tilde{H}_2= -1$ and $\tilde{H}_4\cdot\tilde{H}_4=
-1$.
\end{proof}

\begin{lemma}
The configuration from Figure~\ref{fig:list} with $(k,r)=(1,2)$
together with the requirement that $E_1$ is not a cutting divisor
cannot occur.
\end{lemma}
\begin{proof}
Suppose otherwise.  Assume that $\tilde{H}_1\cdot\tilde{H}_1= -1$ and
$\tilde{H}_3\cdot\tilde{H}_3= -1$.  If this is not the case, then by
Lemmas~\ref{th:atl} and \ref{th:pair} we may assume that
$\tilde{H}_4\cdot\tilde{H}_4= -1$ and $\tilde{H}_2\cdot\tilde{H}_2=
-1$ and argue similarly.  The curves beyond $\tilde{H}_1$ have
self-intersection strictly less than $-1$.  The curve immediately
adjacent and beyond $\tilde{H}_1$ is $\tilde{E}_1$ and this has
self-intersection strictly less than $-2$.  This is because we must
blow up between $E_1$ and $H_4$ to separate cycles, and also between
$\tilde{E}_1$ and $\tilde{L}_1$ to break cycles and to maintain
$\tilde{H}_1\cdot\tilde{H}_1= -1$ and $\tilde{H}_3\cdot\tilde{H}_3=
-1$.  Thus if we blow down $\tilde{H}_1$ the remaining branch beyond
$\tilde{L}_{\infty}$ consists of curves with self-intersection
strictly less than $-1$.  Also $\tilde{H}_2$ has self-intersection
strictly less than $-1$ since we have to blow up the intersection
between $H_2$ and $H_4$ and the intersection between $H_2$ and $L_1$
in order to break cycles and maintain $\tilde{H}_3\cdot\tilde{H}_3=
-1$.  After blowing down $\tilde{H}_1$, $\tilde{L}_{\infty}$ has
self-intersection $0$ and valency $3$ with two branches consisting of
curves of self-intersection strictly less than $-1$.  Thus we can use
Lemma~ \ref{th:neg} to get a contradiction.
\end{proof}

\section{Non-isotrivial rational polynomials of simple type}
\label{sec:class}
The configuration in Figure~\ref{fig:confC} is the starting point for
any non-isotrivial rational polynomial of simple type.  Notice that we
can fill one puncture in each fibre of any such map to get an
isotrivial family of curves and the puncture varies linearly with
$c\in\bbC$.  Notice also that there is an irregular fibre for each of
the $r$ intersection points of the $(1,1)$ curve with $(1,0)$
horizontal curves away from $L_{\infty}$.  In fact there is at most
one more irregular fibre which can only occur in rather special cases,
as we discuss in subsection \ref{subsec:irr}.

From now on the configuration $\C$ is given by Figure~\ref{fig:confC}
with $r+2$ horizontal curves.  Beginning with $\C$ we will list all of
the rational polynomials of simple type generated from this
configuration.  We shall give the splice diagrams for these
polynomials first. Although we compute the polynomials later,
geometric information of interest is often more easily extracted from
the splice diagram or from our construction of the polynomials than
from an actual polynomial.

The splice diagram encodes the topology of the polynomial.  It
represents the link at infinity of the \regular{} fibre, or it can be
thought of as an efficient plumbing graph for the divisor at infinity,
$D\subset X$.  It encodes an entire parametrised family of polynomials
with the same topology of their regular fibres.  See
\cite{ENeThr,NeuCom,NeuIrr} for more details.  Within this family,
polynomials can still differ in the topology of their irregular
fibres.  Our methods also give all information about the irregular
fibres, as we describe in subsection \ref{subsec:irr}.

The configuration $\C$ has $r+3$ irreducible components so when we
blow up to get $D$ by Lemma~\ref{th:propD} ~(i) we will leave $r+1$
exceptional curves behind as cutting divisors.  By
Lemma~\ref{th:propD} ~(ii) we must break the $r$ cycles in $\C$ with
multiple blow-ups at the points of intersection leaving $r$
exceptional curves behind as cutting divisors.  We blow up multiple
times between the $r$th $(1,0)$ horizontal curve and the $(1,1)$
horizontal curve in order to break a cycle.  Thus, we require those
blow-ups to satisfy the condition that the exceptional curve will
break the cycle if removed.  Equivalently, each new blow-up takes
place at the intersection of the most recent exceptional curve with an
adjacent curve.  We call such a multiple blow-up a {\em \toric{}}.

\begin{figure}[ht]
$$\xymatrix@R=12pt@C=15pt@M=0pt@W=0pt@H=0pt
{{\hbox to 0 pt{\hskip2pt P\hss}}&\circ\ar@{--}[dd]\\
&\circ\ar@{--}[r]&&\ar@{--}[l]\circ\ar@{--}[r]\ar@{.}[dd]&&\circ\ar@{--}[l]\\
&\circ&&\ar@{.}[r]&\ar@{--}[r]&\circ\ar@{--}[r]\ar@{--}[d]&&\circ\ar@{--}[l]\\
&&&&{\hbox to 0 pt{\hskip2pt$-1$\hss}}&\circ
}$$
\caption{Sequence of blow-ups starting at $P$ and ending at the $-1$ curve.}
\label{fig:nonbreak}
\end{figure}
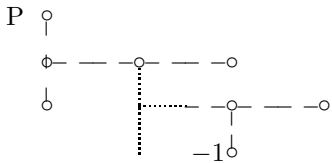

We have one extra cutting divisor. This will arise as the last
exceptional curve blown up in a sequence of blow-ups that does not
break a cycle.  We will call this
sequence of blow-ups a {\em \nonsep{}}. A
priori, this \nonsep{}  could be a
sequence as in Figure~\ref{fig:nonbreak},
where the final $-1$ curve is the cutting divisor. However, we shall
see that the extra nodes this introduces in the dual graph prohibit
$D$ from blowing down to a Morrow configuration, so the sequence is
simply a string of $-2$ exceptional curves followed by $-1$
exceptional curve that is the cutting divisor. This arises from
blowing up a point on a curve in the blow-up of $\C$ that does not lie
on an intersection of irreducible components.  

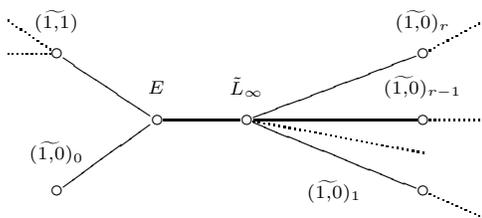
\begin{figure}[ht]  
  $$\xymatrix@R=6pt@C=8pt@M=0pt@W=0pt@H=0pt{
    \ar@{.}[dr]&_{(\widetilde{1,1})}&&&&&&&&
    &_{(\widetilde{1,0})_r}&\ar@{.}[dl]\\
    \ar@{.}[r]&\circ\ar@{-}[ddrrr]&&&&&&&&&\circ\ar@{-}[ddlll]\\
    &&&&_E&&&_{\tilde{L}_{\infty}}&&&_{(\widetilde{1,0})_{r-1}}&
    &&&&&&&\\
    &&&&\circ\ar@{-}[rrr]\ar@{-}[ddlll]&&&\circ\ar@{-}[rrr]
    \ar@{-}[ddrrr]\ar@{.}[drrr]&&&\circ\ar@{.}[r]&\\
    &_{(\widetilde{1,0})_0}&&&&&&&&&&\\
    &\circ&&&&&&&&_{(\widetilde{1,0})_1}&\circ\ar@{.}[dr]&\\
    &&&&&&&&&&&&& }$$
\caption{Dual graph of $\C$ blown up at points of intersection.}
\label{fig:plumb1}
\end{figure}

Let us begin by just performing the \toric s at the points of
intersection, of $\C$ and leaving the \nonsep{}
until later.  This gives the dual graph in Figure~\ref{fig:plumb1}
with the proper transforms of the $r+1$ $(1,0)$ horizontal curves and
the $(1,1)$ horizontal curve indicated along with $\tilde{L}_{\infty}$
and the exceptional curve $E$ arising from the blow-up of the triple
point in $\C$.  There are $r$ branches heading out from the proper
transform of $(1,1)$ consisting of curves of self-intersection less
than $-1$ and beyond each of the proper transforms of the $r$ $(1,0)$
horizontal curves the curves have self-intersection less than $-1$.

The self-intersection of each of $(\widetilde{1,0})_0$, $E$ and
$\tilde{L}_{\infty}$ is $-1$.  The self-intersections of
$(\widetilde{1,1})$ and $(\widetilde{1,0})_i$, $i=1,\dots,r$ are
negative and depend on how we blow up at each point of intersection.

\begin{lemma}   \label{th:oneb}
  There is at most one branch in $D$ beyond $(\widetilde{1,1})$, and
  $r-1$ of the horizontal curves $(\widetilde{1,0})_i$ (those with
  index $i=1,\dots,r-1$ say) have self-intersection $-1$ and only $-2$
  curves beyond.
\end{lemma}
\begin{proof}
  Since the self-intersection of each of the curves beyond
  $(\widetilde{1,1})$ is less than $-1$ each branch beyond
  $(\widetilde{1,1})$ cannot be blown down before $(\widetilde{1,1})$.
  Thus, there are at most two branches.
  
  Furthermore, since the self-intersection of each of the curves
  beyond $(\widetilde{1,0})_i$, $i=1,\dots,r$ is less than $-1$, the
  branch beyond $(\widetilde{1,0})_i$ can be blown down before
  $(\widetilde{1,0})_i$ only if $(\widetilde{1,0})_i$ has
  self-intersection $-1$ and each curve beyond has self-intersection
  $-2$.  Thus, at most two branches beyond $(\widetilde{1,0})_i$,
  $i=1,\dots,r$ do not consist of a $-1$ curve with a string of $-2$
  curves beyond.  If there are two such branches then the blow-ups
  that create them create corresponding branches beyond
  $(\widetilde{1,1})$ (or possibly just decrease the intersection
  number at $(\widetilde{1,1})$). These two branches cannot be fully
  blown down until everything else connecting to the $\tilde L_\infty$
  vertex are blown down, but the vertex $(\widetilde{1,1})$ and any
  branches beyond it cannot blow down first. Thus $D$ cannot blow down
  to a Morrow configuration. Thus there is at most one such branch,
  proving the Lemma.
\end{proof}

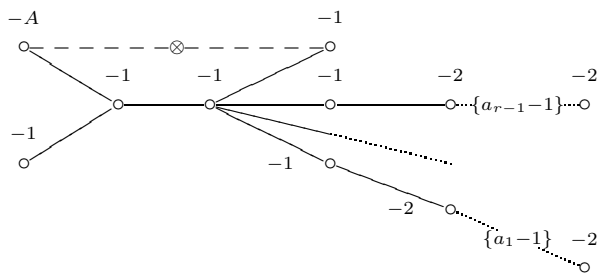
\begin{figure}[ht]  
  $$\xymatrix@R=6pt@C=8pt@M=0pt@W=0pt@H=0pt{
    _{-A}&&&&&&&&&_{-1}\\
    \circ\ar@{-}[ddrrr]\ar@{--}[rrrrrrrrr]|-\otimes&&
    &&&&&&&\circ\ar@{-}[ddlll]\\
    &&&_{-1}&&&_{-1}&&&_{-1}&&&_{-2}&&&&&_{-2}\\
    &&&\circ\ar@{-}[rrr]\ar@{-}[ddlll]&&&\circ\ar@{-}[rrr]
    \ar@{-}[ddrrr]\ar@{-}[drrr]&&&\circ\ar@{-}[rrr]&&&
    \circ\ar@{.}[rrrrr]|-{\{a_{r-1}-1\}}&&&&&\circ\\
    _{-1}&&&&&&&&&\ar@{.}[drrr]\\
    \circ&&&&&&&&_{-1}&\circ\ar@{-}[ddrrr]&&&\\
    \\
    &&&&&&&&&&&_{-2}&\circ\ar@{.}[ddrrrrr]|-{\{a_1-1\}}\\
    &&&&&&&&&&&&&&&&&_{-2}\\
    &&&&&&&&&&&&&&&&&\circ }$$
\caption{Dual graph of partially blown-up configuration of curves.}
\label{fig:plumb2}
\end{figure}

Figure~\ref{fig:plumb2} gives the dual graph of the partially blown up
$\C$ where the label of each curve is now its self-intersection
number.  The branch beyond $(\widetilde{1,0})_i$ consists of a string
of $a_i-1$ $-2$ curves and $A=\sum_{i=1}^{r-1}a_i$.  We have thus far
only blown up once between the $r$th $(1,0)$ horizontal curve and the
$(1,1)$ horizontal curve, indicating the exceptional divisor by
$\otimes$.  We may blow up many more times---perform a 
\toric{}---leaving behind the final exceptional curve as cutting
divisor to get a branch beyond $(\widetilde{1,1})$ and a branch beyond
$(\widetilde{1,0})_r$.  In addition, we still have to perform the
\nonsep{} at some point on the divisor.

\begin{lemma}  \label{th:nonse}
  The \nonsep{} occurs beyond either
  $(\widetilde{1,1})$, $(\widetilde{1,0})_r$, or $(\widetilde{1,0})_0$
  and in the latter case $(\widetilde{1,1})\cdot(\widetilde{1,1})=-1$.
\end{lemma}
\begin{proof}
  If the \nonsep{} occurs on the branch beyond
  $(\widetilde{1,0})_i$, $i=1,\dots,r-1$ then that branch cannot be
  blown down.  By the proof of lemma~\ref{th:oneb}, in order to obtain
  a linear graph we must blow down $r-1$ of the branches beyond
  $(\widetilde{1,0})_i$, $i=1,\dots,r$.  Thus, if the 
  \nonsep{} 
  does occur beyond $(\widetilde{1,0})_i$ for some $i\le
  r-1$, then the $(\widetilde{1,0})_r$ branch blows down, so we simply
  swap the labels $i$ and $r$.
  
  The \nonsep{} cannot occur on $E$ or
  $\tilde{L}_{\infty}$ because the resulting cutting divisor would
  not be sent to a finite value.
  
  If the \nonsep{} occurs on the branch beyond
  $(\widetilde{1,0})_0$ then we must be able to blow down the branch
  beyond $(\widetilde{1,1})$, hence the branch must consist of
  $(\widetilde{1,1})$ with self-intersection $-1$.
\end{proof}

\begin{lemma}  \label{th:crem}
  We may assume the \nonsep{} does not occur beyond
  $(\widetilde{1,0})_0$.
\end{lemma}
\begin{proof}
  By Lemma~\ref{th:nonse} if the \nonsep{} occurs
  beyond $(\widetilde{1,0})_0$ then
  $(\widetilde{1,1})\cdot(\widetilde{1,1})=-1$.  In particular,
  $1=A=\sum_1^{r-1}a_i$.  Thus, $r=2$, $a_1=1$.  With only four
  horizontal curves, we can perform a Cremona transformation to make
  $(\widetilde{1,0})_0$ the $(1,1)$ curve and hence we are in the
  first case of Lemma~\ref{th:nonse}.
\end{proof}

\begin{lemma}  \label{th:nons}
The \nonsep{} occurs on either of the last curves
beyond $(\widetilde{1,1})$ or $(\widetilde{1,0})_r$ and is a string of
$-2$ curves followed by the $-1$ curve that is a cutting divisor.
\end{lemma}
\begin{proof}
  Arguing as previously, if the \nonsep{} occurs
  anywhere else, or if it is more complicated, then it introduces a
  new branch preventing the divisor $D$ from blowing down to a linear
  graph.
\end{proof}

%\begin{lemma}  \label{th:irr}
%There are exactly $r$ irregular fibres except if there is no branch
%beyond $(\widetilde{1,1})$ when there are $r$ or $r+1$ irregular
%fibres.
%\end{lemma}
%\begin{proof}
%The irregular fibres arise from intersection points in $\C$ away from
%$L_{\infty}$ and from the \nonsep{}.  This would be
%$r+1$ irregular fibres except that the \nonsep{}
%usually occurs at (more precisely, infinitely near to) one of the
%intersection points.  Only when there is no branch beyond
%$(\widetilde{1,1})$ can we perform the \nonsep{}
%away from an intersection point and gain another irregular fibre.
%\end{proof}

We now know that our divisor $D$ results from
Figure~\ref{fig:plumb2} by doing a \toric{} between the $(1,1)$
curve and the $r$-th $(1,0)$ curve, leaving behind the final $-1$
exceptional curve as a cutting divisor and then performing a
\nonsep{} on a curve adjacent to this cutting
divisor to produce second cutting divisor. 

A priori, it is not clear that this procedure always gives rise to a
divisor $D\subset X$ where $X$ is a blow-up of $\bbP^2$ and $D$ is the
pre-image of the line at infinity.  The classification will be complete
once we show it does.

\begin{lemma}   \label{th:anyt}
  The above procedure always gives rise to a configuration that blows
  down to a Morrow configuration (see Lemma \ref{th:propD}) and hence
  determines a rational polynomial of simple type.
\end{lemma}
\begin{proof}
  The calculation involves the relation between plumbing graphs and
  splice diagrams described in \cite{ENeThr} or \cite{NeuIrr}, with
  which we assume familiarity. In particular, we use the continued
  fractions of weighted graphs described in \cite{ENeThr}. If one has
  a chain of vertices with weights $-c_0,-c_1,\dots,-c_t$, its
  continued fraction based at the first vertex is defined to be
  $$
  c_0-\cfrac1{c_1-\cfrac 1{c_2-\genfrac{}{}{0pt}{}{}{\ddots~
        \genfrac{}{}{0pt}{}{}{-\genfrac{}{}{0pt}{}{}{\cfrac1{c_t}}}}}}$$

  The dual graph for the curve configuration of Lemma \ref{th:anyt}
  has chains starting at the vertex $(\widetilde{1,1})$ and
  $(\widetilde{1,0})_r$. We claim these chains have continued
  fractions evaluating to $A-1+\frac PQ$ and $\frac qp$ respectively,
  where $P,Q,p,q$ are arbitrary positive integers with $Pq-pQ=1$.  We
  describe the main ingredients of this calculation but leave the
  details to the reader.
  
  An easy induction shows that the initial \toric{}
  leads to chains at $(\widetilde{1,1})$ and $(\widetilde{1,0})_r$
  with continued fractions $A-1+\frac nm$ and $\frac mn$ with positive
  coprime $n$ and $m$. The \nonsep{} then changes the
  fraction $\frac nm$ or $\frac mn$ that it operates on as follows. If
  the \nonsep{} consists of $k$ blow-ups at the end
  of the left chain then $\frac nm$ is replaced by $\frac NM$ with
%  $Nm-nM=1$ and $kn<N\le (k+1)n$. If the \nonsep{}
  $Nm-nM=1$ and $k\le \frac Mm<\frac Nn\le (k+1)$. If the
  \nonsep{} is on the right then $\frac mn$ is
  similarly changed instead.
  
  Renaming, we can describe this in terms of our chosen names
  $p,q,P,Q$ as follows. We either have $P>p$ or $q>Q$. If $P>p$ the
  initial \toric{} leads to chains with continued
  fractions $A-1+\frac pq$ and $\frac qp$ and the 
  \nonsep{}
  then consists of a sequence of $k:=\lfloor\frac Qq\rfloor$
  blowups extending the left chain (and changing its continued
  fraction to $A-1+\frac PQ$). If $q>Q$ the continued fractions are
  $A-1+\frac PQ$ and $\frac QP$ after the separating blowup and the
  non-separating blow-up consists of $k:=\lfloor\frac pP\rfloor$
  blow-ups extending the right chain (and changing its continued
  fraction to $\frac qp$).
  
  To prove the Lemma we must show that the dual graph of our curve
  configuration blows down to a Morrow configuration. We can blow down
  the chains starting at $(\widetilde{1,0})_i$, $i=0,\dots,r-1$, to
  get a chain. To check that this chain is a Morrow configuration we
  must compute its determinant, which we can do with continued
  fractions as in \cite{ENeThr}. We first replace the two end chains
  by vertices with the rational weights $-A+1-\frac PQ$ and $-\frac
  qp$ determined by their continued fractions to get a chain of four
  vertices with weights
  $$
  -A+1-\frac PQ,\quad 0,\quad -1+A,\quad -\frac qp.$$
  Then,
  computing the continued fraction for this chain based at its right
  vertex gives $\frac qp-\frac QP = \frac{Pq-pQ}{Pp}=\frac1{Pp}$,
  showing that the determinant is $-1$ as desired, and completing the
  proof.
\end{proof}
%% END OF COMMENTING OUT

\begin{thm}\label{th:main}
  Given positive integers $P,Q,p,q$ with $Pq-pQ=1$ and positive
  integers $a_1,\dots,a_{r-1}$, the splice diagram of our rational
  polynomial $f$ of simple type with non-isotrivial fibres is given in
  Figure~\ref{fig:splice1} with
  \begin{align*}
  A&=a_1+\dots+a_{r-1},\\
  B&=AQ+P-Q,\\
C&=Aq+p-q,\\
b_i&=qQa_i+1 \quad\text{for each $i$.}
  \end{align*}
The degree of $f$ is:\quad $\operatorname{deg}(f)=A(Q+q)+P+p$.
\end{thm}

\def\Dot{\lower.2pc\hbox to 2pt{\hss$\bullet$\hss}}
\def\Circ{\lower.2pc\hbox to 2pt{\hss$\circ$\hss}}
\def\Vdots{\raise5pt\hbox{$\vdots$}}
\begin{figure}[ht]
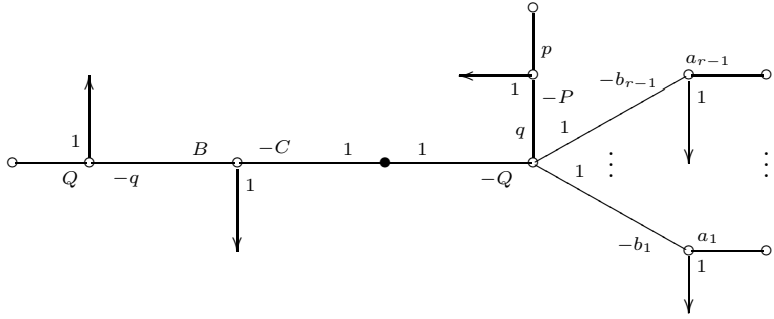

\iffalse
$$
\objectmargin{0pt}\spreaddiagramrows{-3pt}
\spreaddiagramcolumns{3pt}\diagram
&&&&&&&&&\Circ\rline^(.25){a_1}\dto^(.25)1&\Circ\\
\Circ\rline^(.75)Q &\Circ\dto^(.25)1\rrline^(.25){-q}^(.75)B&&
\Circ\dto^(.25)1\rrline^(.25){-C}^(.75)1&&\Dot\rrline^(.25)1^(.75){-Q}
&&\Circ\dline_(.25)q_(.75){-P}\urrline^(.25)1^(.75){-b_1}
\drrline^(.25)1_(.75){-b_r}&\Vdots&&\Vdots\\
&&&&&&&\Circ\dline^(.25)p\lto^(.25)1&&\Circ\rline^(.25){a_r}
\dto^(.25)1&\Circ\\
&&&&&&&\Circ&&&
\enddiagram
$$
\fi
$$
\objectmargin{0pt}\spreaddiagramrows{-3pt}
\spreaddiagramcolumns{3pt}\diagram
&&&&&&&\Circ&&&\\
&&&&&&&\Circ\uline_(.35)p\lto^(.25)1&&\Circ\rline^(.25){a_{r-1}}\dto^(.25)1&\Circ\\
\Circ\rline_(.75)Q &\Circ\uto^(.25)1\rrline_(.25){-q}^(.75)B&&
\Circ\dto^(.25)1\rrline^(.25){-C}^(.75)1&&\Dot\rrline^(.25)1_(.75){-Q}
&&\Circ\uline^(.35)q_(.75){-P}\urrline^(.25)1^(.75){-b_{r-1}}
\drrline^(.25)1_(.75){-b_1}&\Vdots&&\Vdots\\
&&&&&&&&&\Circ\rline^(.25){a_1}
\dto^(.25)1&\Circ\\
&&&&&&&&&&
\enddiagram
$$
\caption{Splice diagram for non-isotrivial rational polynomial.}
\label{fig:splice1}
\end{figure}
(In \cite{NNoMo} an ``additional'' case was given, which is, however,
of the above type with $P=Q=p=1$, $q=a_1=2$.)
\begin{proof}
  For the following computations we continue to assume the reader is
  familiar with the relationship between resolution graphs and splice
  diagrams described in \cite{NeuIrr}.  The arrows signify places at
  infinity of the \regular{} fibre, one on each horizontal curve.  The
  fact that $(\widetilde{1,0})_r$ is next to $\tilde{L}_{\infty}$ in
  the dual graph says that the edge determinant of the intervening
  edge is $1$.  This corresponds to the fact that $Pq-pQ=1$, which we
  already know.  Similarly, $(\widetilde{1,0})_i$ is next to
  $\tilde{L}_{\infty}$ for $i=1,\dots,r-1$ so the weight $b_i$ is
  determined by the edge determinant condition $b_i=qQa_i+1$.
  The
  ``total linking number'' at the vertex corresponding to each
  horizontal curve (before blowing down $(\widetilde{1,0})_0$) is zero
  (terminology of \cite{NeuIrr}; this reflects the fact that the link
  component corresponding to the horizontal curve has zero linking
  number with the entire link at infinity, since at almost all points
  on a horizontal curve, the polynomial has no pole). The weight $C$
  is determined by the zero total linking number of
  $(\widetilde{1,1})$, giving $C=Aq+p-q$. For any $i$ the fact that
  vertex $(\widetilde{1,0})_i$ has zero total linking gives
  $B=AQ+P-Q$.
\end{proof}

It is worth summarising some consequences of our construction that
will be useful later.
\begin{lemma}
  The number of blow-ups in the final \nonsep{} is
  %$k:=\max(\lceil\frac Pp\rceil,\lceil\frac qQ\rceil)-1$ and these
  $k:=\max(\lfloor\frac Qq\rfloor,\lfloor\frac pP\rfloor)$ and these
  blow-ups occurred at the $(Q,-q)$ branch or the $(p,-P)$ branch of
  the above splice diagram according as the first or second entry of
  this {\rm max} is the larger. Moreover, the non-separating blow-ups
  occurred on the corresponding horizontal curve if and only if $q=1$
  resp.\ $P=1$.
\end{lemma}
\begin{proof}
  The first part was part of the proof of Lemma \ref{th:anyt}.  For
  the second part, note that if $q=1$ then certainly $q>Q$ must fail,
  so $P>p$ and the nonseparating blow-ups were on the left. The
  continued fraction on the left was $A-1+\frac pq=A-1+p$ which is
  integral, showing that the left chain consisted only of the
  exceptional curve before the non-separating blow-up. Conversely, if
  the non-separating blow-ups were adjacent to that exceptional curve
  then the left chain was a single vertex, hence had integral
  continued fraction, so $q=1$. The argument for $P=1$ is the same.
\end{proof}
%The polynomial has no multiple fibres in the special cases of
%$(P,Q,p,q)=(2,2b-1,1,b)$ or $(P,Q,p,q)=(1,b,1,b+1)$ and in these
%cases, if each $a_i>1$ then the example is generically non-singular.

\begin{thm}\label{th:defspace}
  The moduli space of polynomials $f\colon\bbC^2\to\bbC$ with the
  above regular splice diagram, modulo left-right equivalence (that
  is, the action of polynomial automorphisms of both domain $\bbC^2$
  and range $\bbC$), has dimension $r+k-2$ with $k$ determined in the
  previous Lemma. In fact it is a $\bbC^k$-fibration over the
  $(r-2)$-dimensional configuration space of $r-1$ distinct points in
  $\bbC^*$ labelled $a_1,\dots,a_{r-1}$, modulo permutations that
  preserve the labelling and transformations of the form $z\mapsto
  az$. 
\end{thm}
\begin{proof}
  The splice diagram prescribes the number of horizontal curves and
  the \toric s at each point of intersection.  The
  only freedom is in the placement of the horizontal curves in
  $\bbP^1\times\bbP^1$, and in the choice of points, on prescribed
  curves, on which to perform the string of blow-ups we call the
  \nonsep{}. The $(1,1)$ horizontal curve is
  \emph{a priori} the graph of a linear map $y=ax+b$ but can be
  positioned as the graph of $y=x$ by by an automorphisms of the image
  $\bbC$.
  
  The point in the configuration space of the Theorem determines the
  placement of the horizontal curves $(1,0)_1,\dots,(1,0)_{r}$ (after
  putting the $(1,0)_0$ curve at $\bbP^1\times\{\infty\}$ and the
  $(1,0)_r$ curve at $\bbP^1\times\{0\}$).  The fibre $\bbC^k$
  determines the sequence of points for the 
  \nonsep{}.
  
  This proves the Theorem, except that we need to be careful, since
  some diagrams occur in the form of Theorem \ref{th:main} in two
  different ways, which might seem to lead to disconnected moduli
  space. But the only cases that appear twice have four horizontal
  curves and the configurations $\C$ are related by Cremona
  transformations.
\end{proof}
This completes the classification of non-isotrivial rational
polynomials of simple type. 

\subsection{The irregular fibres}\label{subsec:irr}

We can read off the topology of the irregular fibres of the polynomial
$f$ of Theorem \ref{th:main} from our construction, since any such
fibre is the proper transform of a vertical $(0,1)$ curve together
with any exceptional curves left behind as cutting divisors when
blowing up on this vertical curve.

We shall use the notation $\bbC(r)$ to mean $\bbC$ with $r$ punctures
(so $\bbC^*=\bbC(1)$), and for the purpose of this subsection we used
$C\cup C'$ to mean disjoint union of curves $C$ and $C'$, and $C+C'$
to mean union with a single normal crossing. The \regular{} fibre of $f$
is $\bbC(r+1)$.

The irregular fibres of $f$ arise through the breaking of cycles
between the $(1,1)$ curve and the $(1,0)_i$ curve for $i=1,\dots,r$,
so there are $r$ of them. The non-separating blow-up also
contributes, but it usually contributes to the $r$-th irregular fibre.
However, if $P=1$ or $q=1$ then the non-separating blow-up occurs on a
horizontal curve and can thus have any $f$-value, so it generically
leads to an additional $(r+1)$-st irregular fibre.

The irregular fibres are all reduced except for the $r$-th irregular
fibre, which is always non-reduced unless one of $P,Q,p,q$ is $1$.

We first assume $q\ne1$ and $P\ne 1$, so there are exactly $r$
irregular fibres.  Then for each $i=1,\dots,r-1$ the $i$-th irregular
fibre is $\bbC(r-1)+\bbC^*$ if $a_i=1$ and $\bbC(r)\cup\bbC^*$ if
$a_i>1$.  The $r$-th irregular fibre is $\bbC(r)\cup\bbC^*\cup \bbC$
generically. As mentioned above, this fibre is reduced if and only if
$Q=1$ or $p=1$. There is a single parameter value in the $\bbC^k$
factor of the parameter space of Theorem \ref{th:defspace} for which
the $r$-th irregular fibre has different topology, namely
$\bbC(r)\cup(\bbC+\bbC)$. In this case it is non-reduced even if $Q=1$
or $p=1$.

If $q=1$ or $P=1$ then write $\frac PQ$ and $\frac qp$ as $\frac 1a$
and $\frac{ak+1}k$ in some order. The non-separating blow-up creates
irregularity in a fibre which generically is distinct from the the
first $r$ irregular fibres. The generic situation is that the $r$-th
irregular fibre is $\bbC(r)\cup\bbC^*$ or $\bbC(r-1)+\bbC^*$ according
as $a>1$ or $a=1$ and the $(r+1)$-st irregular fibre is
$\bbC(r+1)\cup\bbC$ or $\bbC(r)+\bbC$ according as $k>1$ or $k=1$, and
both are reduced. But there are codimension $1$ subspaces of the
parameter space for which the topology is different. For instance, the
$(r+1)$-st irregular fibre will be non-reduced if one blows up more
than once on a vertical curve while doing the 
\nonsep{} that
creates it.

\subsection{Monodromy} \label{subsec:monodromy} 

We can also read off the monodromy for our polynomial $f$.  Consider a
generic vertical $(0,1)$ curve $C$ in our construction.  Removing its
intersections with the horizontal curves gives a regular fibre $F$ of
$f$. Since we have positioned the horizontal curve $(1,0)_0$ at
$\infty$ we think of $F$ as an $r+1$-punctured $\bbC$. We call the
intersection of the $(1,1)$ horizontal curve with $C$ the $0$-th
puncture of $F$ and for $i=1,\dots,r$ we call the intersection of the
$(1,0)_i$ curve with $C$ the $i$-th puncture of $F$.

If the $(r+1)$-st irregular fibre exists the local monodromy around it
is trivial.  For $i=1,\dots,r$ the monodromy around the $i$-th
irregular fibre rotates the $0$-th puncture of the regular fibre
$\bbC(r+1)$ around the $i$-th puncture.  In terms of the braid group
on the $r+1$ punctures, with standard generators $\sigma_i$ exchanging
the $(i-1)$-st and $i$-th puncture for $i=1,\dots,r$, the local
monodromies are $h_1=\sigma_1^2$,
$h_2=\sigma_1\sigma_2^2\sigma_1^{-1}$, \dots,
$h_r=\sigma_1\dots\sigma_{r-1}\sigma_r^2\sigma_{r-1}^{-1}
\dots\sigma_1^{-1}$.  The monodromy $h_\infty=h_r\dots h_1$ at
infinity is $\sigma_1\sigma_2\dots\sigma_r\sigma_r\dots\sigma_1$.
It is not hard to verify that $h_1,\dots,h_r$ freely generate a free
subgroup of the braid group.

\section{Explicit polynomials}\label{sec:poly}
The splice diagram gives sufficient information (Newton polygon,
topological properties, etc.) that one can easily find the polynomial
without significant computation by making an educated guess and then
confirming that the guess is correct.  The answer is as follows:

{\bf Case 1.} $k\le\frac pP<k+1$.
(Then $\frac pP<\frac qQ\le k+1$.)

Let $s_1=\alpha_0+\alpha_1x+\dots+\alpha_{k-1}x^{k-1}+x^ky$. Let
$\beta_1,\dots\beta_{r-1}$ be distinct complex numbers in $\bbC^*$.

$$f(x,y)=x^{q-Qk}s_1^Q+x^{p-Pk}s_1^P
\prod_{i=1}^{r-1}(\beta_i-x^{q-Qk}s_1^Q)^{a_i}.$$

{\bf Case 2.} $k\le\frac Qq<k+1$. (Then $\frac Qq<\frac Pp\le k+1$.)

Let $s_2=\alpha_0+\alpha_1y+\dots+\alpha_{k-1}y^{k-1}+xy^k$. Let
$\beta_1,\dots\beta_{r-1}$ be distinct complex numbers in $\bbC^*$.

$$f(x,y)=y^{Q-qk}s_2^q+y^{P-pk}s_2^p
\prod_{i=1}^{r-1}(\beta_i-y^{Q-qk}s_2^q)^{a_i}.$$

One can compute the splice diagram and see it is correct. One can
verify that the \regular{} fibres are rational by the explicit
isomorphism:
$$f^{-1}(t) \to \bbC-\{0,\beta_1,\dots,\beta_{r-1}, t\}\qquad\left\{
\begin{matrix}
(x,y)&\mapsto &x^{q-Qk}s_1^Q&&&\text{(Case 1),}\\ 
(x,y)&\mapsto &y^{Q-qk}s_2^q&&&\text{(Case 2),}  
\end{matrix}\right.$$
for generic
$t$. The irregular values of $t$ are
$0,\beta_1,\dots,\beta_{r-1}$ if $P\ne1$ and $q\ne1$. If $P=1$ then
$t=\alpha_0\prod \beta_i$ is the additional irregular value that our
earlier discussion predicts, and if $q=1$ then $t=\alpha_0$ is the
additional irregular value.

The space of parameters
$(\alpha_0,\dots,\alpha_{k-1},\beta_1,\dots,\beta_{r-1})$ maps to the
moduli space we computed earlier with fibre of dimension $1$.  Indeed,
with $B,C$ as in Theorem \ref{th:main}, the polynomial
$$f_\lambda(x,y)=\lambda^{-1}f(\lambda^Bx,\lambda^{-C}y)$$ has the same
form with the parameters $\beta_j$ replaced by $\lambda^{-1}\beta_j$
and $\alpha_j$ replaced by $\lambda^{jB+A-1}\alpha_j$.

\vspace{6pt}
To put the above polynomials in the form of $f_1(x,y)$ of Theorem
\ref{th:summary}, in case 1 we 
rename the exponents
$q-Qk$ to $q_1$, $p-Pk$ to $p_1$, $Q$ to $q$, $P$ to $p$.
In case 2 we rename $Q-qk$ to $q_1$, $P-pk$ to $p_1$, and then
exchange $x$ and $y$.

\section{The isotrivial case.}\label{sec:miyanishi-sugie}
After the first version of this paper was completed we realised that
the classification in \cite{MSuGen} for the isotrivial case has
omissions. In this section we therefore sketch the corrected
classification using the techniques of this paper. The discussion of
the parameter spaces and the irregular fibres for the resulting
polynomials is similar to the non-isotrivial case, so we leave it to
the reader. One can give an alternative proof using Kaliman's
classification \cite{KalPol} of all isotrivial polynomials.

We will restrict ourselves to the case of ample rational polynomials,
i.e. those with at least three $(1,0)$ horizontal curves.  The case of
one $(1,0)$ horizontal curve always gives a polynomial equivalent to a
coordinate by the Abhyankar-Moh-Suzuki theorem \cite{AMoEmb,SuzPro}.
The case of two $(1,0)$ horizontal curves is dealt with from a splice
diagram perspective in \cite{NeuCom} and earlier by analytic methods
in \cite{Saito}. The result is included in our summary Theorem
\ref{th:summary}.

As before, compactify $\bbC^2$ to $X$ and construct a map
$X\rightarrow\bbP^1\times\bbP^1$.  The map is essentially canonical
(up to an automorphism of one factor.)  The image of the divisor at
infinity $D\subset X$ in $\bbP^1\times\bbP^1$ is given by a collection
of $(1,0)$ curves since we used three of the horizontal curves to get
a map to $\bbP^1\times\bbP^1$ and in order that the fibres give an
isotrivial family, any other horizontal curves must also be $(1,0)$
curves.

When there are at least three $(1,0)$ horizontal curves, by the
following lemma the original configuration of curves in
$\bbP^1\times\bbP^1$ breaks into the two cases of no vertical curves
or one vertical curve.

\begin{lemma}
  An ample rational polynomial with isotrivial fibres has at most one
  vertical curve over a finite value.
\end{lemma}
\begin{proof}
  We can argue as in the previous section.  The curve over infinity,
  $L_{\infty}$ is not blown up since there are no triple points.  If
  there is more than one vertical curve over a finite value then there
  are precisely three $(1,0)$ horizontal curves since otherwise there
  would be at least two $(1,0)$ horizontal curves that would be blown
  up at least twice and since all curves beyond these horizontal
  curves (exceptional curves or vertical curves) have
  self-intersection $<-1$ we would get two branches $B_1$ and $B_2$
  made up of the proper transforms of these two $(1,0)$ horizontal
  curves and all curves beyond these, meeting at a valency $>2$ curve,
  $L_{\infty}$, with self-intersection $0$.  This is the impossible
  situation of Lemma~\ref{th:neg}.

There can be at most two vertical curves since if there are $l$
vertical curves we need to break $2l$ cycles but since there are
precisely three $(1,0)$ horizontal curves, we begin with $l+4$ curves
so we can break at most $l+2$ cycles by Lemma~\ref{th:propD} ~(i).
Therefore $2l\leq l+2$ so $l\leq 2$.

The lemma follows when we get rid of the case of two vertical curves
and three $(1,0)$ horizontal curves.  The few cases are easily
dismissed by hand.
\end{proof}

So the beginning configuration is given by Figure~\ref{fig:conf7} or
Figure~\ref{fig:conf8}. We analyse these below as Case 1 and Case 2.
\begin{figure}[ht]
$$\xymatrix@R=12pt@C=15pt@M=0pt@W=0pt@H=0pt{
\ar@{-}&\ar@{-}[ddddddd]&&&&&&\\
\ar@{-}[rrrrrr]&&&&&&&\\
\ar@{-}[rrrrrr]&&&&&&&\\
\ar@{-}[rrrrrr]&&&&&&&\\
\ar@{-}[rrrrrr]&&&\ar@{}[dd]^(.3){.}^{.}^(.7){.}&&&&\\
\ar@{-}&&&&&&&\\
\ar@{-}[rrrrrr]&&&&&&&\\
&}$$
\caption{Configuration of horizontal curves with $L_{\infty}$.}
\label{fig:conf7}
\end{figure}
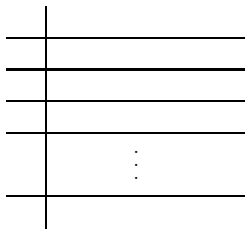
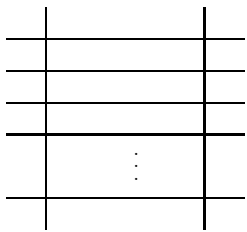
\begin{figure}[ht]
$$\xymatrix@R=12pt@C=15pt@M=0pt@W=0pt@H=0pt{
\ar@{-}&\ar@{-}[ddddddd]&&&&\ar@{-}[ddddddd]&&\\
\ar@{-}[rrrrrr]&&&&&&&\\
\ar@{-}[rrrrrr]&&&&&&&\\
\ar@{-}[rrrrrr]&&&&&&&\\
\ar@{-}[rrrrrr]&&&\ar@{}[dd]^(.3){.}^{.}^(.7){.}&&&&\\
\ar@{-}&&&&&&&\\
\ar@{-}[rrrrrr]&&&&&&&\\
&&&&&&&}$$
\caption{Configuration of horizontal curves with $L_{\infty}$ and a 
vertical curve over a finite value.}
\label{fig:conf8}
\end{figure}

\subsection*{Case 1.}  Denote by $r$ the number of horizontal
curves.
In Figure~\ref{fig:conf7} we must leave behind $r-1$
curves as cutting divisors. To do so we do a 
\nonsep{}
on each of $r-1$ horizontal curves (anything else leads to a
configuration of curves whose intersection matrix has determinant $0$,
and which can therefore not blow down to a Morrow configuration).
Thus, on the $i$-th horizontal curve we blow up $a_i$ times and then
leave behind the final exceptional divisor, giving a string of $-2$
curves of length $a_i-1$.

\begin{figure}[ht]
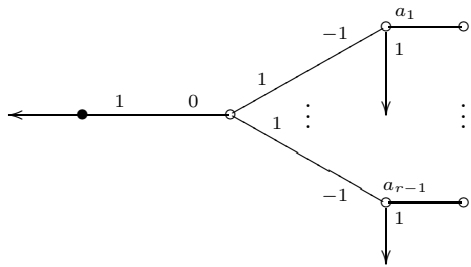

$$
\objectmargin{0pt}\spreaddiagramrows{-3pt}
\spreaddiagramcolumns{3pt}\diagram
&&&&&\Circ\rline^(.25){a_1}\dto^(.25)1&\Circ\\
&\Dot\lto\rrline^(.25)1^(.75){0}
&&\Circ\urrline^(.25)1^(.75){-1}
\drrline^(.25)1_(.75){-1}&\Vdots&&\Vdots\\
&&&&&\Circ\rline^(.25){a_{r-1}}
\dto^(.25)1&\Circ\\
&&&&&&
\enddiagram
$$
\caption{Splice diagram for case 1 of isotrivial fibres.}
\label{fig:splice3}
\end{figure}

The resulting splice diagram is as in Figure~\ref{fig:splice3}.

This splice diagram has been analysed in \cite{NeuIrr}, where it is shown
that its general polynomial is
$$f(x,y)=y\prod_{i=1}^{r-1}(x-\beta_i)^{a_i}+h(x),$$
where $h(x)$ is a polynomial of degree $<\sum_{i=1}^{r-1}a_i$.

This case covers the following cases from \cite{MSuGen}: Case 1 of
Theorem 3.3., Theorem 3.7, Case I of Theorem 3.10.

\subsection*{Case 2.}
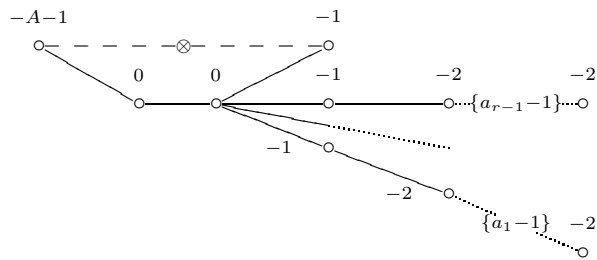
\begin{figure}[ht]  
  $$\xymatrix@R=6pt@C=8pt@M=0pt@W=0pt@H=0pt{
    _{-A-1}&&&&&&&&&_{-1}\\
    \circ\ar@{-}[ddrrr]\ar@{--}[rrrrrrrrr]|-\otimes&&
    &&&&&&&\circ\ar@{-}[ddlll]\\
    &&&_{0}&&&_{0}&&&_{-1}&&&_{-2}&&&&&_{-2}\\
    &&&\circ\ar@{-}[rrr]&&&\circ\ar@{-}[rrr]
    \ar@{-}[ddrrr]\ar@{-}[drrr]&&&\circ\ar@{-}[rrr]&&&
    \circ\ar@{.}[rrrrr]|-{\{a_{r-1}-1\}}&&&&&\circ\\
    &&&&&&&&&\ar@{.}[drrr]\\
    &&&&&&&&_{-1}&\circ\ar@{-}[ddrrr]&&&\\
    \\
    &&&&&&&&&&&_{-2}&\circ\ar@{.}[ddrrrrr]|-{\{a_1-1\}}\\
    &&&&&&&&&&&&&&&&&_{-2}\\
    &&&&&&&&&&&&&&&&&\circ }$$
\caption{Dual graph of partially blown-up configuration of curves for 
Fig.~\ref{fig:conf8}.}
\label{fig:ms1}
\end{figure}

Denote by $r+1$ the number of horizontal curves.
In Figure~\ref{fig:conf8} we must do \toric s at
$r$ intersection points and then do an additional 
\nonsep{}. 
As in the Section \ref{sec:class}, one finds that each of
$r-1$ of the \toric s creates a string of $-2$ curves
attached to the corresponding horizontal curve, while the last one can
be arbitrary, as described in the proof of Lemma
\ref{th:anyt}. In Figure~\ref{fig:ms1} we show the situation after
doing the first $r-1$ \toric s and doing the first
step of the $r$-th one.

Moreover, the \nonsep{} then occurs adjacent to the
exceptional curve left behind in the final \toric{}.
The analysis is almost identical to the proof of Lemma \ref{th:anyt},
with the resulting strings now having continued fractions $A+\frac PQ$
and $\frac qp$ respectively, with notation as in that proof.

The resulting splice diagram is as in Figure~\ref{fig:ms2}, with
notation exactly as in Theorem \ref{th:main}.
\def\Vdots{\raise5pt\hbox{$\vdots$}}
\begin{figure}[ht]
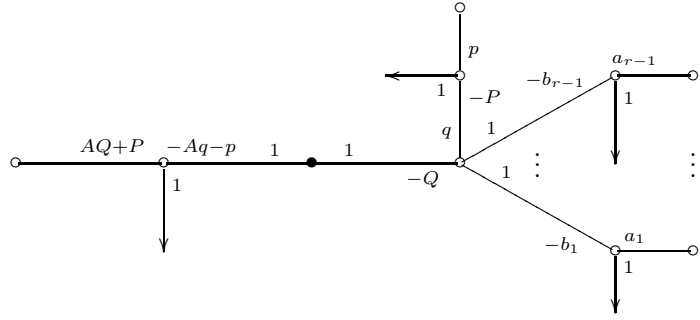

$$
\objectmargin{0pt}\spreaddiagramrows{-3pt}
\spreaddiagramcolumns{3pt}\diagram
&&&&&&\Circ&&&\\
&&&&&&\Circ\uline_(.35)p\lto^(.25)1&&\Circ\rline^(.25){a_{r-1}}\dto^(.25)1&\Circ\\
\Circ\rrline^(.65){AQ+P}&&
\Circ\dto^(.25)1\rrline^(.25){-Aq-p}^(.75)1&&\Dot\rrline^(.25)1_(.75){-Q}
&&\Circ\uline^(.35)q_(.75){-P}\urrline^(.25)1^(.75){-b_{r-1}}
\drrline^(.25)1_(.75){-b_1}&\Vdots&&\Vdots\\
&&&&&&&&\Circ\rline^(.25){a_1}
\dto^(.25)1&\Circ\\
&&&&&&&&&
\enddiagram
$$
\caption{Splice diagram for Case 2 of isotrivial fibres.}
\label{fig:ms2}
\end{figure}
The polynomial in this case is exactly as in Section \ref{sec:poly}
except that the first term $x^{q-Qk}s_1^Q$ respectively
$y^{Q-qk}s_2^q$ is omitted. Namely, let
$\beta_1,\dots\beta_{r-1}$ be distinct complex numbers in $\bbC^*$ and
let $k$ be as in Theorem \ref{th:defspace}. 

If $k\le\frac pP<k+1$
(so $\frac pP<\frac qQ\le k+1$),
let $s_1=\alpha_0+\alpha_1x+\dots+\alpha_{k-1}x^{k-1}+x^ky$. Then
$$f(x,y)=x^{p-Pk}s_1^P
\prod_{i=1}^{r-1}(\beta_i-x^{q-Qk}s_1^Q)^{a_i}.$$

If $k\le\frac Qq<k+1$ (so $\frac Qq<\frac Pp\le k+1$),
let $s_2=\alpha_0+\alpha_1y+\dots+\alpha_{k-1}y^{k-1}+xy^k$. Then
$$f(x,y)=y^{P-pk}s_2^p
\prod_{i=1}^{r-1}(\beta_i-y^{Q-qk}s_2^q)^{a_i}.$$

This case covers the following cases from \cite{MSuGen}: Cases 2,3,4
of Theorem 3.3 and Case II of Theorem 3.10.  However, \cite{MSuGen}
only has examples in which one of $P,Q,p,q$ is equal to $1$.
 
Note that the isotrivial splice diagrams of Case 1
and Case 2 can be considered to belong to one family: putting
$(P,Q)=(1,0)$ in Figure~\ref{fig:ms2} gives Figure~\ref{fig:splice3}.
Nevertheless, the two cases have rather different geometric
properties.

\section{General rational polynomials.}
In this section we will give a result for ample rational polynomials
that are not necessarily of simple type.

\begin{prop}  \label{th:hor-1}
An ample rational polynomial contains a $(1,0)$ horizontal curve whose
proper transform has self-intersection $-1$ and meets
$\tilde{L}_{\infty}$.
\end{prop}
\begin{proof} 
By the classification of ample rational polynomials of simple type,
the proposition is true in this case.  So, we may assume that there is
a horizontal curve of type $(m,n)$ for $m>1$.

Suppose there is no $(1,0)$ horizontal curve with the property of the
proposition.  Then by the proof of Lemma~\ref{th:triple} there are at least two
$(1,0)$ horizontal curves whose proper transforms have
self-intersection $<-1$ and meet $\tilde{L}_{\infty}$.  
%As in the proof of Lemma~\ref{th:inunion} 
By Lemma~\ref{le:beyond} any curves beyond these horizontal curves
have self-intersection $<-1$.

A horizontal curve of type $(m,n)$ must meet $L_{\infty}$ at exactly
one point, and hence with a tangency of order $m$ or at a singularity
of the curve.  This is because if a horizontal curve were to meet
$L_{\infty}$ twice then we would not be able to break cycles since
when we blow up next to $L_{\infty}$, those exceptional curves are
sent to infinity under the polynomial and hence must be retained
in the configuration of curves.  Thus we must blow up there to
get a configuration of curves with normal intersections.  The final
exceptional curve in such a sequence of blow-ups will have
self-intersection $-1$ and valency $>2$.

If we can blow down the configuration of curves then eventually at
least one curve adjacent to the $-1$ curve is blown down and hence the
$-1$ curve ends up with non-negative self-intersection.  But the final
configuration is not a linear graph since the proper transforms of the
two $(1,0)$ horizontal curves and any curves beyond give two branches.
Thus the final configuration is not a Morrow configuration which
contradicts Lemma \ref{th:propD}.
\end{proof}

The following result is a generalisation of Lemma~\ref{th:n=1}.
\begin{cor}
For any ample rational polynomial, a smooth horizontal curve of type
$(m,n)$ with $m>0$ must be of type $(m,1)$.
\end{cor}
\begin{proof}  
The statement is true for $m=1$ by Lemma~\ref{th:n=1} so will assume
$m>1$.  A curve of type $(m,n)$ will intersect the $(1,0)$ horizontal
curves $m$ times, with multiplicity, unless possibly if the $(m,n)$
curve is singular at these points of intersection.  The latter
possibility is ruled out by the assumption of the corollary.  Hence
the $(1,0)$ horizontal curves will be blown up at least $m$ times and
their proper transforms will have self-intersection $<-m$.  This
contradicts the previous proposition so the result follows.
\end{proof}

When the rational polynomial is not ample, Russell has an example of a
horizontal curve of type $(3,2)$.  See the examples in the next
section.  Note that smoothness of the horizontal curve is necessary in
the corollary (at the points of intersection with the $(1,0)$
horizontal curves) since we can always have two horizontal curves of
types $(l,1)$ and $(m,1)$ and together they can be considered as a
singular horizontal curve of type $(l+m,2)$.

\subsection{Adding horizontal curves}

Consider the following construction on $\bbC^2$. Blow up repeatedly
starting at a point on the $y$-axis so that the resulting exceptional
curves form a chain from the $y$-axis to the last exceptional curve
blown up.  If we now remove the $y$-axis and all but the last
exceptional curve from the blown-up $\bbC^2$ we get a new $\bbC^2$
that we call $\bbC^2_1$.  Any polynomial $f\colon\bbC^2\to\bbC$
induces a polynomial $f_1\colon\bbC^2_1\to\bbC$.  Suppose the $y$-axis
intersects generic fibres of $f$ in $d$ points.  Then the \regular{}
fibres of $f_1$ are simply \regular{} fibres of $f$ with $d$ extra
punctures. In fact, this construction simply adds an extra degree $d$
horizontal curve, namely the $y$-axis becomes a degree $d$ horizontal
curve for $f_1$.

From the point of view of the polynomials, what we have done is
replaced $f(x,y)$ by 
$$f_1(x,y)=f(x,s),\quad s=a_0+a_1x+\dots+a_{k-1}x^{k-1}+x^ky,$$
that is, we have composed $f$ with the birational morphism
$(x,y)\mapsto (x,s)$ of $\bbC^2$. 

Since one can compose $f$ first with a polynomial automorphism to
raise its degree, one can easily add horizontal curves of arbitrarily high
degree by this construction. This makes clear that any classification
of non-simple-type polynomials must take account of this sort of
operation, including composition with more complicated birational
morphisms.

Although this is a complication, it can also simplify some issues.

Here is a simple illustrative example. We start with the simplest
rational polynomial $g(x,y)=x$, apply a polynomial automorphism to get
$f(x,y)=x+y^2$ and then apply the above birational morphism to get
$f_1(x,y)=x+(a_0+a_1x+\dots+a_{k-1}x^{k-1}+x^ky)^2$ with one degree
one horizontal and one degree two horizontal. It is not hard to check
(e.g., by listing possible splice diagrams) that this gives, up to
equivalence, the only non-simple-type polynomials with \regular{} fibre
$\bbC-\{0,1\}$, so with the classification of simple type polynomials,
we get:
\begin{prop}
  A polynomial with general fibre $\bbC-\{0,1\}$ is left-right
  equivalent to one of the form
  $f_2(x,y)$ or $f_3(x,y)$ of Theorem 
  \ref{th:summary} with $r=2$ or $r=3$ respectively, or to 
$f(x,y)=x+(a_0+a_1x+\dots+a_{k-1}x^{k-1}+x^ky)^2$.\qed
\end{prop}

This proposition also follows from 
Kaliman's classification \cite{KalPol} of isotrivial
polynomials.

\section{Examples}
It is worth including some interesting known examples of rational
polynomials from the perspective used in this paper.  These examples
are neither of simple type nor ample.

\begin{figure}[ht]
\begin{center}
\scalebox{0.6}{\includegraphics{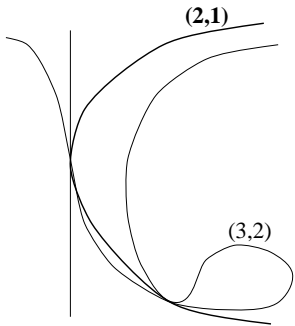}}
\end{center}
\caption{A bad field generator.}
\label{fig:rus}
\end{figure}

Russell \cite{RusGoo} (correctly presented in \cite{BCaOne})
constructed an example of a rational polynomial with no degree one
horizontal curves.  This is an example of a bad field generator---a
polynomial that is one coordinate in a birational transformation but
not in a birational morphism.  It is given by beginning with three
curves in $\bbP^1\times\bbP^1$ as in Figure~\ref{fig:rus}.  The
$(2,1)$ curve and the $(3,2)$ curve intersect at an order three
tangency and at the same point the $(3,2)$ intersects itself at a
tangency.  They are the two horizontal curves of the polynomial.  The
vertical curve is $L_{\infty}$.

The actual polynomial in this case is, with $s=xy+1$, $$
f(x,y)=(y^2s^4+y(s+xy)s+1)(ys^5+2xys^2+x)$$
and the splice diagram is 
$$\xymatrix@R=12pt@C=12pt@M=0pt@W=0pt@H=0pt{
  \Circ\ar@{-}[rr]^(.75){3}&&\Circ\ar[ld]\ar[d]\ar[rd]
  \ar@{-}[rrr]^(.25){-4}&&&
   \Dot\ar@{-}[rrr]^(.75){-2}&&&\Circ\ar@{-}[rr]^(.25){3}
\ar@{-}[dd]^(.75){-13}
  &&\Circ\\ &&&&&&&&&&\\
&&&&&&\Circ\ar[lu]\ar[ld]\ar@{-}[rr]^(.25){-27}&&
\Circ\ar@{-}[dd]^(.25){2}\\
&&&&&&&&&&\\
&&&&&&&&\Circ}$$

\begin{figure}[ht]
\begin{center}
\scalebox{0.6}{\includegraphics{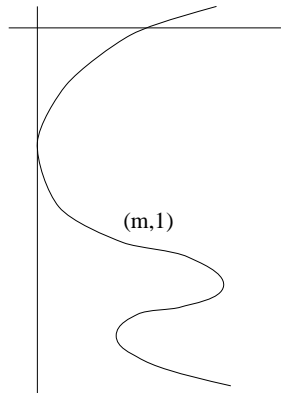}}
\end{center}
\caption{Classification of rational polynomials with a $\bbC^*$ fibre.}
\label{fig:kal}
\end{figure}

Kaliman \cite{KalRat} classified all rational polynomials with one
fibre isomorphic to $\bbC^*$.  Figure~\ref{fig:kal} gives three curves
in $\bbP^1\times\bbP^1$, the two horizontal curves and $L_{\infty}$.
The $(m,1)$ curve has the property that when it is mapped downwards
onto a $(1,0)$ curve, there are only two points of ramification, both
with maximal ramification of $m$, at $L_{\infty}$ and at the irregular
fibre isomorphic to $\bbC^*$.  Kaliman's entire classification begins
with this configuration of curves.  The only points that can be blown
up are those that are infinitely near to the point of intersection of
the two horizontal curves (besides the unnecessary blowing up where
the $(m,1)$ curve meets $L_{\infty}$) and one exceptional curve is
left behind as a component of the reducible fibre.

\vfill

\iffalse
Recall from Section~\ref{sec:class} that the classification of
polynomials with \regular{} fibre isomorphic to an annulus lies behind
the classification of rational polynomials of simple type with
non-isotrivial fibres.  It seems reasonable to conjecture that the
classification of polynomials with \regular{} fibre isomorphic to an
annulus along with Kaliman's classification of rational polynomials
with one fibre isomorphic to $\bbC^*$ can be used in a similar way to
classify all ample rational polynomials with one horizontal curve of
degree greater than one.
\fi

%\clearpage

\end{document}